\documentclass[11pt,leqno,a4paper]{article}

\usepackage{xy} 
\xyoption{all} 
\newdir^{ (}{{}*!/-3pt/\dir^{(}}     
\newdir^{  }{{}*!/-3pt/\dir^{}}     
\newdir_{ (}{{}*!/-3pt/\dir_{(}}     
\newdir_{  }{{}*!/-3pt/\dir_{}}

\def\theenumi{(\alph{enumi})}

\def\p@enumii{\theenumi} 
 
\usepackage{amsmath} 
 
\usepackage{amssymb} 
 
\usepackage{amscd} 
 
\usepackage{ifthen}

\newcounter{zahl} 
 
\newenvironment{aufz}{\begin{list}{{\rm \alph{zahl})}} {\usecounter{zahl} \labelwidth0.8cm \labelsep0.2cm \leftmargin1.0cm \topsep0.0cm \itemsep0.1cm \parsep0.0cm}}{\end{list}}

\newcommand{\DS}{\displaystyle} 
 
\newcommand{\TS}{\textstyle} 
 
\newcommand{\SC}{\scriptstyle}

\DeclareMathOperator{\Hom}{Hom} 
 
\DeclareMathOperator{\Frob}{Frob}

\DeclareMathOperator{\End}{End} 
 
\DeclareMathOperator{\Isom}{Isom} 
 
\DeclareMathOperator{\Aut}{Aut} 
 
\DeclareMathOperator{\im}{im}

\DeclareMathOperator{\CKoh}{\check H}

\DeclareMathOperator{\Cov}{Cov}

\DeclareMathOperator{\Spec}{Spec} 
 
\DeclareMathOperator{\BSpec}{{\bf Spec}} 
 
\DeclareMathOperator{\Spf}{Spf} 
 
\DeclareMathOperator{\Spm}{Sp} 
 
\DeclareMathOperator{\GL}{GL}

\DeclareMathOperator{\Lie}{Lie} 
 
\DeclareMathOperator{\Stab}{Stab}

\DeclareMathOperator{\id}{\,id} 
 
\DeclareMathOperator{\Id}{Id}

\newcommand{\rig}{{\rm rig}}

\newcommand{\alg}{{\rm alg}}

\newcommand{\topol}{{\rm top}} 
 
\newcommand{\temp}{{\rm temp}}

\newcommand{\tors}{{\rm tors}}

\newcommand{\et}{{\rm \acute{e}t}}

\DeclareMathOperator{\rk}{rk}

\DeclareMathOperator{\charakt}{char}

\DeclareMathOperator{\Res}{Res}

\renewcommand{\phi}{\varphi} 
\renewcommand{\epsilon}{\varepsilon}

\usepackage{amsfonts} 
 
\newcommand{\BOne}{\underline{\hbox{\rm1\kern-2.5pt l\kern.9pt}}} 
\newcommand{\OOne}{\BO\kern-7.7pt\raisebox{1.3pt}{\underline{\phantom{{$\BO$}}}}}

\newcommand{\Cz}{\mathbb{C}}

\newcommand{\Fz}{\mathbb{F}} 
 
\newcommand{\Gr}{\mathbb{G}}

\newcommand{\Nz}{\mathbb{N}} 
 
\renewcommand{\Pr}{\mathbb{P}}

\newcommand{\Rz}{\mathbb{R}}

\newcommand{\rbij}{\mbox{\mathsurround=0pt \;$-\hspace{-0.75em}\stackrel{\sim\quad}{\longrightarrow}$\;}}

\newcommand{\es}{\enspace}

\newcommand{\open}{^\circ} 
 
\newcommand{\mal}{^{^\times}}

\newcommand{\dpl}{\mbox{\rm (\hspace{-0.13em}(}} 
 
\newcommand{\dpr}{\mbox{\rm )\hspace{-0.13em})}} 
 
\newcommand{\invlim}[1][]{\ifthenelse{\equal{#1}{}}
{\DS \lim_{\longleftarrow}}
{\DS \lim_{\underset{#1}{\longleftarrow}}}
} 
 
\newcommand{\dirlim}[1][]{\ifthenelse{\equal{#1}{}}
{\DS \lim_{\longrightarrow}}
{\DS \lim_{\underset{#1}{\longrightarrow}}}
} 
 
\newcommand{\ul}[1]{{\underline{#1}}} 
 
\newcommand{\nul}[1]{{\underline{#1}}} 
 
\newcommand{\ol}[1]{{\overline{#1}}}

\newcommand{\wh}[1]{{\widehat{#1}}} 
 
\newcommand{\wt}[1]{{\widetilde{#1}}}

\usepackage{amsthm} 
 
\theoremstyle{plain} 
 
\newtheorem{thm}{Theorem}[section]

\newtheorem{lemma}[thm]{Lemma} 
 
\newtheorem{keylemma}[thm]{Key-Lemma}

\newtheorem{cor}[thm]{Corollary} 
 
\newtheorem{prop}[thm]{Proposition}

\theoremstyle{definition} 
 
\newtheorem{Def}[thm]{Definition}

\newtheorem{Bigexample}[thm]{Example} 
 
\newtheorem{question}[thm]{Question} 
 
\newtheorem{remark}[thm]{Remark} 
 
\newtheorem{example}[thm]{Example} 
 
\newtheorem{examples}[thm]{Examples} 
 
\def\?{\ ???\ \immediate\write16{}%
 
\immediate\write16{Warning: There was still a question mark . . . }%
 
\immediate\write16{}}

\long\def\forget#1{}

\def\CA{\mathcal{A}} 
 
\def\CE{\mathcal{E}} 
\def\CF{\mathcal{F}} 
\def\CG{\mathcal{G}} 
\def\CH{\mathcal{H}} 
 
\def\CM{\mathcal{M}} 
\def\CO{\mathcal{O}} 
\def\CS{\mathcal{S}} 
\def\CV{\mathcal{V}} 
\def\FM{\mathfrak{M}} 
 
\def\BO{\mathbb{O}} 
\def\BC{\mathbb{C}} 
\def\BF{\mathbb{F}} 
\def\BL{\mathbb{L}} 
\def\BH{\mathbb{H}} 
\def\BG{\mathbb{G}} 
 
\def\Fa{\mathfrak{a}} 
\def\Fp{\mathfrak{p}} 
\def\Fm{\mathfrak{m}} 
\def\Fn{\mathfrak{n}} 
\def\FA{\mathfrak{A}} 
\def\FC{\mathfrak{C}} 
 
\def\BN{\mathbb{N}} 
\def\BZ{\mathbb{Z}} 
\DeclareMathOperator{\Max}{Max} 
\DeclareMathOperator{\Tor}{Tor} 
\DeclareMathOperator{\Cent}{Cent} 
 
\DeclareMathOperator{\Sym}{Sym} 
 
\def\Ltt{{L\langle t\rangle}} 
\newcommand{\UCF}{\underline{\CF}\mathstrut} 
\newcommand{\UCE}{\underline{\CE}\mathstrut} 
\newcommand{\UCG}{\underline{\CG}\mathstrut} 
\newcommand{\UCH}{\underline{\CH}\mathstrut}

\def\longto{\longrightarrow} 
\def\into{\hookrightarrow} 
 
\newcommand{\dash}{{\textrm{-}}\,} 
\DeclareMathOperator{\bVec}{{\bf Vec}} 
 
\title{Uniformizable families of $t$-motives}

\begin{document}

\author{\footnotesize by\\ \\ 
\hbox to 4cm{\hfill Gebhard B{\"o}ckle}\ \  
{\footnotesize and}\ \  
\hbox to 4cm{Urs Hartl\hfill} 
} 
 
\maketitle 
 
 
\begin{abstract} 
\noindent 
Abelian $t$-modules and the dual notion of $t$-motives were introduced by
Anderson as a generalization of Drinfeld modules. For such Anderson defined
and studied the important concept of uniformizability. It is an interesting question, and
the main objective of the present article to see how uniformizability behaves
in families. Since uniformizability is an analytic notion, we have to work
with families over a rigid analytic base. 
We provide many basic results, and in fact a large part of
this article concentrates on laying foundations for studying the above
question. Building on these, we obtain a generalization of a uniformizability
criterion of Anderson and, among other things, we establish that the locus of
uniformizability is Berkovich open.

\end{abstract} 
 
\section*{Introduction} 
 
Over the complex numbers $\BC$, any elliptic curve $E$ is biholomorphically isomorphic 
to a quotient $\BC/\Lambda$, where $\Lambda$ is the period lattice of $E$. The analogue of the complex numbers 
in characteristic $p$ is the field $\BC_\infty$ which is the topological closure of an algebraic closure 
of $\BF_q\dpl 1/t\dpr$. With this analogy in mind, Drinfeld associated to any discrete finitely generated  
$\BF_q[t]$-sublattice $\Lambda$ of $\BC_\infty$ an algebraic object, which we now call a Drinfeld-module. 
 
In the case of elliptic curves, the uniformization $\BC/\Lambda$ is not just a point-wise property. 
For any complex analytic family $E\to S$ of elliptic curves over a complex manifold $S$, if $E_0$ over a point 
$s_0$ is given by a lattice $\Lambda_0$, then this construction may be extended locally analytically 
to a lattice over an open neighborhood of $s_0$ in~$S$. In fact, there is a tautological analytical local system 
$\BL:=(\Lambda_z\subset\BC)_{z\in\BH}$ on the upper half plane $\BH$ which represents a uniformization of the  
universal family of elliptic curves on $\BH$. Now for any simply connected neighborhood $U\subset S$ of $s_0$, 
there exists an analytic morphism $U\to\BH$ such that $E$ restricted 
to $U$ is obtained from the universal family on $\BH$ via pullback, and thus the 
corresponding local system is the pullback of~$\BL$. 
The same can be done for Drinfeld-modules of rank $r$ using the rigid analytic space $\Omega^r$ instead of $\BH$, 
and where one considers the rigid analytic Grothendieck topology instead of the complex analytic one, cf.~\cite{drinfeld}  
and \cite{boeckle}, Chap.~4. 
 
\smallskip 
 
In \cite{anderson}, Anderson introduced the notion of a $t$-module, generalizing that 
of a Drinfeld-module. Again there arises the question of the uniformizability of 
these objects, i.e., the question of whether one can always write these in a suitable sense as 
$\BC_\infty^d/\Lambda$ for some discrete lattice $\Lambda$. As Anderson shows, not 
all $t$-modules can be written in such a way. If they can be, he calls them \emph{uniformizable}, and  
moreover any uniformizable $t$-module has good reduction. 
 
When compared with the classical situation, the nature of uniformizability cannot be explained  
satisfactorily. One thinks of $t$-modules as being analogous to abelian varieties.  But 
over $\BC$ any abelian variety is isomorphic to $\BC^g$ modulo a suitable discrete lattice of  
full rank~$2g$. Over a $p$-adic complete algebraically closed field $\BC_p$ only abelian  
varieties with bad, but totally multiplicative reduction admit a uniformization by a discrete  
lattice. 
 
\medskip 
 
To any $t$-module, Anderson associates a dual object called $t$-motive, and the 
uniformizability of the $t$-module is equivalent to the analytic triviality of the associated $t$-motive. 
We explain the correspondence in the last section of the article. 
Moreover $t$-motives are a special case of so-called $\tau$-sheaves. 
The aim of this article is to study analytic triviality in families of $\tau$-sheaves. 
 
The first example of a family of $\tau$-sheaves arises from a family of Drinfeld modules on a $\BC_\infty$-variety $X$. Since all Drinfeld modules are uniformizable, every fiber of this family is analytically trivial. Moreover, as we have remarked above the family of Drinfeld modules can even be uniformized as a whole on a covering of $X$. However, this covering will in general not be a variety but a rigid analytic space. 
 
Beyond Drinfeld modules, the next interesting examples were constructed by 
Gardeyn in his thesis, cf.~\cite[II.2]{GardeynDiss}.  
In Section~\ref{SectExample} of this article, we describe 
another non-trivial family which was computed by R.\ Pink and which in fact describes a  
fine moduli space of pure ``polarized'' $t$-motives of rank and dimension~$2$. In these families not every fiber is analytically trivial. 
 
Given a family of $\tau$-sheaves on a variety $X$ over $\BC_\infty$, the ultimate goal is to describe the nature of 
the locus of analytic triviality inside $X$. We expect that the equations (better: inequalities) describing analytic triviality are rarely algebraic -- and in the example of Section~\ref{SectExample} this is indeed not the case. Therefore to investigate analytic triviality 
we are again urged to work in a rigid analytic setting. 
In this article we show that the locus of analytic triviality is Berkovich  
open, cf.~Theorem~\ref{UnivOpen}. However the known examples suggest that in fact the locus of analytic triviality  
of a family of $\tau$-sheaves on an affinoid space is the complement of a quasi-compact rigid analytic subspace. This is a much stronger assertion.  
 
\smallskip 
 
To fix ideas, let $R$ denote an affinoid algebra over $\BC_\infty$,  
and let $(M,\tau)$ be a pair of a finitely generated projective module $M$ over the Tate-algebra  
$R\langle t\rangle$ over $R$, and an endomorphism $\tau:M\to M$ such that  
$\tau( rt^i m)=t^ir^q\tau(m)$ for $r\in R$ and $m\in M$. (This is the  
affine version of a rigid analytic $\tau$-sheaf, cf.~Definition~\ref{DefRigTau}.) The module $M$ gives  
rise to a quasi-coherent sheaf $\CF$ on $X:=\Spm R$, and $\tau$ is an $\BF_q[t]$-linear operation on it. 
Let $W(\CF)$ denote the associated {\'e}tale sheaf on the {\'e}tale site over~$X$, and consider 
the left exact \emph{fundamental sequence} 
\begin{equation}\label{ASlikeSeq}0\longto \UCF^\tau\longto 
  W(\CF)\stackrel{1-\tau}\longto W(\CF), 
\end{equation} 
which defines the {\'e}tale sheaf $\UCF^\tau$. The latter is a sheaf of $\BF_q[t]$-modules. 
Its construction is reminiscent of Artin-Schreier theory (and specializes to it, for instance  
for $t\mapsto 0$). The pair $(M,\tau)$ is called \emph{analytically trivial}  
if {\'e}tale locally the rank of $\UCF^\tau$  
as a free module over $\BF_q[t]$ is the same as the rank of $M$ over~$R\langle t\rangle$. 
 
It turns out that $\UCF^\tau$ is overconvergent 
in the sense of \cite{JongPut}. Therefore, we work on the rigid analytic site not 
just with the usual classical points, but with so-called \emph{analytic points}, first defined 
by van der Put and later developed by Berkovich. These points can detect exactness in sequences of  
overconvergent sheaves. Some necessary definitions on this are recalled in an appendix. 
In Section~\ref{RigTauSec} we define our basic objects of study -- rigid analytic $\tau$-sheaves. 
Equipped with these tools and definitions, we shall prove the following results: 
 
\smallskip

In Section~\ref{KeyLemSec} we establish a key lemma on the exact sequence (\ref{ASlikeSeq}).  
Its first consequence is the overconvergence of $\UCF^\tau$. 
It also implies that the sequence is right exact if and only if it is so at all
{\'e}tale analytic stalks. Moreover the $\BF_q[t]$-modules of sections of $\UCF^\tau$ are 
free of rank at most that of $M$, and for $V\to U$ an {\'e}tale cover with $V$ connected, the inclusion 
$\UCF^\tau(U)\into\UCF^\tau(V)$ is injective and saturated. 
 
Then in Section~\ref{SectTriv}, we investigate triviality of $\tau$-sheaves, where the pair $(M,\tau)$ is trivial, 
if $M$ is generated over $R\langle t\rangle$ by the global sections $\UCF^\tau(X)$. For  
$R$ a field our results are translations of Anderson from \cite{anderson}. The case of a point is the  
crucial step when we prove a simple condition for global triviality.  
Analytic triviality is clearly a weaker condition than triviality, and we study this 
concept in detail in the subsequent Section~\ref{SectAnalytTriv}. There we prove that 
the analytically trivial locus is Berkovich open. We also show that on the one hand analytic triviality is  
equivalent to analytic triviality at all analytic stalks, and on the other that it is equivalent to  
the triviality of the $\tau$-sheaf over a suitable temperate (in the sense of \cite{andre2}) {\'e}tale covering. 
 
In Section~\ref{GenArtSch} we investigate natural conditions under which the sequence~(\ref{ASlikeSeq}) is right exact. Essentially we 
require that the $\tau$-sheaf is point-wise composed of trivial and nilpotent $\tau$-sheaves. Section~\ref{SectExample} 
is dedicated to an example by R.\ Pink in which he computes the analytically trivial locus of a $2$-dimensional moduli space 
of $t$-modules which are not Drinfeld-modules. 
 
The last section, Section~\ref{SectFamAMod}, describes the transition from $t$-modules to $\tau$-sheaves,  
describes the resulting dualities on the level of lattices, and formulates and proves Anderson's  
uniformizability criterion for families of $t$-motives. Here we also prove the existence of Anderson's exponential 
map over an arbitrary (rigid analytic) base.

\medskip 
 
{\bf Acknowledgements:} We would both like to thank R.\ Pink for his interest
in this work, his example and many related mathematical discussions. Our
thanks also goes to the anonymous referee for some insightful comments that
led to a simplification in the proof of Proposition~\ref{ExampleProp}.

 
\section{Rigid Analytic $\tau$-Sheaves}\label{RigTauSec} 
\setcounter{equation}{0} 
 
In this section, we introduce our basic objects of study, namely rigid 
analytic families of $\tau$-sheaves. They were first introduced in 
\cite{boeckle} and form a natural generalization of analytic objects 
defined and studied first by Anderson in \cite[{\S}~2]{anderson}. 
We recall those definitions and constructions which will be relevant  
to the question of analytic triviality, which is 
the central topic of this article.

 
\subsection{Definitions} 
 
Let $\Fz_q$ be the field of $q$ elements and characteristic $p$, let 
$C$ be a complete, smooth, geometrically irreducible curve over 
$\Fz_q$, and $\infty$ a fixed closed point on $C$. Let $A$ be the ring 
of regular functions on $C\smallsetminus \{\infty\}$ and $K$ its field 
of fractions. The completion of $K$ at the place $\infty$ is denoted 
$K_\infty$, and the completion of an algebraic closure of $K_\infty$ 
is denoted $\Cz_\infty$. For an $\Fz_q$-scheme $X$ denote by 
$\sigma_X$ its absolute Frobenius endomorphism with respect to $\BF_q$  
which acts as the identity on points and as the $q$-power map on the  
structure sheaf.  
 
\smallskip 
 
In the following definition, the ring $\tilde A$ is either $A$ or 
$A/\Fa$ for some non-zero ideal $\Fa$ of~$A$. 
\begin{Def} 
For $X$ as above an {\em (algebraic) $\tau$-sheaf over $\tilde A$ on $X$\/} is a pair $\UCF:=(\CF,\tau_\CF)$ consisting of a coherent sheaf $\CF$ on $X \times_{\Fz_q} \Spec \tilde{A}$ and an $\CO_{X \times_{\Fz_q} \Spec \tilde{A}}$-module homomorphism 
\[ 
\tau_\CF:\enspace(\sigma_X \times \id)^\ast \,\CF \longrightarrow \CF\,. 
\] 
\end{Def} 
 
 
\smallskip 
 
The notion of $\tau$-sheaf as above was first introduced in 
\cite{boeckle-pink}. An earlier version that required local freeness 
was defined in \cite{taguchi-wan}. In  \cite{boeckle-pink} it 
is also shown that the category of $\tau$-sheaves is stable under  
tensor product, and pullback and proper pushforward along morphisms 
in~$X$. 
 
\smallskip 
 
As the present article centers around the notion of 
uniformizability, we need, as in \cite{boeckle}, to pass from the 
above algebraic to a corresponding rigid analytic setting. (See Bosch, G{\"u}ntzer, Remmert \cite{BGR} for a general introduction to rigid analytic geometry.) 
Let $L$ be a complete non-archimedean valued field with $K_\infty 
\subset L\subset \BC_\infty$. If $X$ is a scheme locally of finite 
type over $L$, we denote by $X^{\rig}$ the associated rigid analytic 
space, cf.~
\cite[{\S} 9.3.4]{BGR}. If $X$ is a 
scheme locally of finite type over $\BF_q$, then we often simply write 
$X^{\rig}$ for $(X\otimes_{\BF_q}L)^{\rig}$. Corresponding to 
$X\mapsto X^\rig$ one has 
a functor $\CF\mapsto\CF^\rig$ on coherent sheaves. If one applies 
rigidification to an algebraic $\tau$-sheaf as above, it is a simple 
matter to verify that one will obtain what we are about to define, 
namely a rigid analytic $\tau$-sheaf. (In fact one obtains two, 
one for each of the coefficient rings we will introduce.) 
 

\smallskip 
 
Unless indicated otherwise, from now on by $X$ we denote a rigid analytic space on~$L$ and by $\sigma_X$ the Frobenius on $X$ with respect to $\BF_q$, i.e., $\sigma_X$ is the $q$-power map on affinoid algebras. We consider the following two natural coefficient rings: 
\begin{enumerate} 
\item $A(\infty):=\Gamma(\FA(\infty), \CO_{\FA(\infty)})$ denotes the ring of entire functions on $\mathfrak{A}(\infty):=(\Spec A)^{\rig}$. 
\item  $A(1):=\Gamma(\FA(1), \CO_{\FA(1)})$ denotes the ring of  rigid analytic functions on the ``unit disc'' $\mathfrak{A}(1)$ of  $\FA(\infty)$, which is constructed as follows:   
\end{enumerate} 
Let $\CO_L$ be the valuation ring of $L$ and consider the formal completion of the scheme $\Spec A \otimes_{\Fz_q} \CO_L$ along its special fiber over the maximal ideal of $\CO_L$. This is an admissible formal scheme over $\Spf \CO_L$ in the sense of Bosch, L{\"u}tkebohmert \cite{FRG}. Then $\mathfrak{A}(1)$ is the rigid analytic $L$-space associated to this formal scheme.  
 
Equivalently, one may also fix a monomorphism $\Fz_q[t] \to A$, denote by $L \langle t\rangle$ the Tate-algebra over $L$ in the variable $t$, and set $A(1) := A \otimes_{\Fz_q[t]} L\langle t\rangle$ and $\mathfrak{A}(1) = \Spm A(1)$. 
 
\medskip 
 
In the following $\FC$ will stand for either $\FA(\infty)$ or $\FA(1)$ 
and $\CA$ for $\Gamma(\FC,\CO_\FC)$. By $\sigma_{L,\mathfrak{C}}$ we 
denote the pullback of the Frobenius $\sigma_{\Spm L}$ on $L$ along 
$\mathfrak{C} \to \Spm L$, i.e., $\sigma_{L,\mathfrak{C}}$ acts as the 
$q$-power map on the coefficients in $L$ and as the identity on the 
variables of $A$. On $X \times_L \mathfrak{C}$ we define the 
endomorphism $\sigma_{X,\mathfrak{C}} := \sigma_X \times 
\sigma_{L,\mathfrak{C}}$. It replaces $\sigma_X \times \id$ in the 
rigid analytic setting. 
 
\begin{Def} \label{DefRigTau} 
Let $X$ be a rigid analytic space over $L$. 
A {\em rigid analytic $\tau$-sheaf over $\CA$ on $X$\/} is a pair $\UCF=(\CF,\tau_\CF)$ consisting of a rigid analytic coherent sheaf $\CF$ of $\CO_{X \times_L \mathfrak{C}}$-modules on $X \times_L \mathfrak{C}$ and an $\CO_{X \times_L \mathfrak{C}}$-module homomorphism 
\[ 
\tau_\CF:\enspace\sigma_{X,\mathfrak{C}}^\ast \,\CF \longrightarrow \CF\,. 
\] 
 
A {\em homomorphism\/} of rigid analytic $\tau$-sheaves over $\CA$ on $X$ 
is a homomorphism of sheaves on $X \times_L \mathfrak{C}$ which commutes 
with the action of $\tau$. It is a {\em monomorphism, epimorphism, 
  isomorphism}, respectively, if its underlying homomorphism of sheaves on 
$X \times_L \mathfrak{C}$ has this property. 
 
A rigid analytic $\tau$-sheaf $\UCF$ is called {\em nilpotent}, it there 
exists an $n\in\BN$ such that $\tau^n\!:(\sigma_{X,\mathfrak{C}}^n)^\ast 
\,\CF \longrightarrow \CF$ is the zero homomorphism. 
\end{Def} 
 
\begin{examples}\label{FirstExs} 
\begin{aufz} 
\item 
Let $\BOne_{X,\CA}$ denote the rigid analytic $\tau$-sheaf over $\CA$ 
on $X$ consisting of the structure sheaf $\CO_{X \times_L 
  \mathfrak{C}}$ and the natural isomorphism  
\[ 
\sigma_{X,\mathfrak{C}}^\ast \,\CO_{X \times_L \mathfrak{C}} \rbij \CO_{X \times_L \mathfrak{C}} \,. 
\] 
 
\item 
Let $\Omega^r := \Pr^r(\Cz_\infty) \smallsetminus  \{\,K_\infty\text{-rational hyper planes}\,\}$ be the {\em $r$-dimensional 
  Drinfeld upper half space\/} and $E$ the universal Drinfeld-module 
on $\Omega^r$. To $E$ is associated an algebraic $\tau$-sheaf $(M(E),\tau)$ 
over $A$, namely the $A$-motive of $E$. The assignment is explained in 
Section~\ref{SectFamAMod}, and further details are given in \cite{boeckle}, Chaps.~4 and~7. 
\end{aufz} 
\end{examples} 
 
\medskip 
 
A rigid analytic $\tau$-sheaf $\UCF$ is called {\em locally free} if 
the underlying sheaf $\CF$ is locally free on $X \times_L 
\mathfrak{C}$, i.e., if there is an admissible covering $U_i$, $i\in 
I$, of $X \times_L \mathfrak{C}$ which trivializes $\CF$. Both of the 
above two examples define in fact locally free $\tau$-sheaves.  
If further $X$ is connected, we denote the {\em rank of $\CF$} by $\rk 
\CF$.  
 
One has the following important result of L{\"u}tkebohmert~\cite[Satz 1]{lubo}: 
\begin{lemma}\label{lubolemma} 
If $\CF$ is a locally free sheaf on $X\times_L \Spm \Ltt$, then there exists an admissible affinoid covering $U_i$ of $X$ such that $\CF|_{U_i \times \Spm\Ltt}$ is free on $U_i\times_L\Spm\Ltt$. 
\end{lemma} 
 
Using the above lemma, one may introduce local coordinates on locally 
free rigid analytic $\tau$-sheaves as follows: 
\begin{remark}[On local coordinates]\label{Taguchi's Trick} 
Let $\UCF$ be a locally free rigid analytic $\tau$-sheaf over 
$A(1)$ on $X$, and fix any non-constant homomorphism $\BF_q[t]\to A$ (it will automatically 
be flat and finite). This induces a finite homomorphism $L\langle t\rangle \to A(1)$, making $A(1)$ into a free $L\langle t\rangle$-module. Thereby $\CF$ becomes a locally free sheaf on $X \times_L \Spm L\langle t\rangle$, say of rank $r$. By Lemma~\ref{lubolemma}, there exists an admissible affinoid covering of $X$ that trivializes $\CF$, i.e., on each affinoid $\Spm B\subset X$ of this covering $\CF$ is associated to a free module $M \cong B\langle t\rangle^r$. On the affinoid algebra $B$ we fix once and for all a complete $L$-algebra norm $|\;\;|$. For any matrix $\Delta$ with entries in $B$ we let $|\Delta|$ be the maximum over the norms of the entries. Choosing a basis $\{\,f_1, \ldots, f_r\,\}$ we write the elements of $M$ as vectors 
\[ 
\sum_{n =0}^\infty b_n t^n\,,\qquad b_n \in B^r\, \quad \text{with} \enspace |b_n| \to 0\enspace \text{for} \enspace n\to \infty\,. 
\] 
We denote by $\sigma$ the action on $B\langle t\rangle$ given by $\sigma\bigl(\sum b_n t^n\bigr) = \sum (b_n)^qt^n$. Then the homomorphism $\tau_\CF$ can be represented with respect to our basis by a matrix 
\[ 
\tau_\CF = \Delta\cdot \sigma\,, \qquad\Delta = \sum_{n=0}^\infty \Delta_n t^n \enspace \in M_r(B\langle t\rangle)\,. 
\] 
Hence $\Delta_n \in M_r(B)$ with $|\Delta_n| \to 0$ for $n\to \infty$. If we replace the basis vectors $f_i$ by $\wt{f}_i = \alpha \, f_i$ for $\alpha \in L\mal$, the $b_n$ are replaced by $\wt{b}_n = \alpha^{-1} \,b_n$ and therefore the matrices $\Delta_n$ are replaced by $\wt{\Delta}_n = \alpha^{1-q}\, \Delta_n$. Hence, we can adjust $|\Delta_n|$. 
 
If there is an element $a\in A\smallsetminus \Fz_q^{alg}$ such that 
$\tau_\CF$ is bijective on $\CF/a\CF$, and if we take $\Fz_q[t] \to 
A:t\mapsto a$ for the above homomorphism, then $\Delta_0$ will be in $\GL_r(B)$. 
Thus if we perform the finite {\'e}tale base extension given by 
adjoining the solution $D$ of the equation $D = \Delta_0 
\,{}^{\sigma\!} D$ and take the columns of $D$ as a new basis, we may further assume that $\Delta_0 =\Id_r$. 
\end{remark}

 
\subsection{Functorial Constructions} 
\label{Functors}

{\em Change of coefficients I\/}: 
 
Let $B$ be an affinoid $L$-algebra and consider $X=\Spec B$ as a 
scheme over $\Fz_q$. Let $\UCF$ be an algebraic $\tau$-sheaf 
over $A$ on $X$. To it one naturally assigns a rigid analytic $\tau$-sheaf 
$(\CF^\rig,\tau_{\CF^\rig})$ over $A(\infty)$ on $\Spm B$ as follows: The 
underlying coherent rigid analytic sheaf $\CF^\rig$ is the pullback of 
$\CF$ along the analytification morphism 
\[ 
X\times_L\mathfrak{C} = (X\times_{\Fz_q} \Spec A)^{\rig} \longrightarrow X \times_{\Fz_q} \Spec A. 
\] 
On it $\tau_\CF$ induces the homomorphism  
\[ 
\tau_{\CF^\rig}\!:\es 
\sigma_{X,\mathfrak{C}}^\ast \,\CF^\rig \enspace = \enspace  \bigl( (\sigma \times \id)^\ast \CF \bigr)^\rig \enspace \xrightarrow{\;\tau_\CF^\rig\;} \enspace \CF^\rig\,. 
\] 
This construction generalizes in an obvious way to arbitrary rigid 
analytic spaces $X$. 
 
\bigskip 
 
\noindent{\em Change of coefficients II\/}: 
 
Let $\UCF$ be a rigid analytic $\tau$-sheaf over 
$A(\infty)$ on $X$. Then the action of $\tau_\CF$ preserves the restriction 
$\CF|_{X\times_L\FA(1)}$, and so we obtain a $\tau$-sheaf over 
$A(1)$ on $X$, which we denote by $\UCF\otimes_{A(\infty)}A(1)$. 
 
\bigskip 
 
\noindent{\em Change of coefficients III\/}: 
 
Let $\UCF$ be a rigid analytic $\tau$-sheaf over 
$A(1)$ on $X$. For any ideal $\Fa$ of $A$ the ideal $\Fa 
A(1)\subset A(1)$ is closed in the topology and satisfies 
$A/\Fa\cong A(1)/\Fa A(1)$. This identification induces an 
obvious functor from rigid analytic $\tau$-sheaves $\UCF$ over 
$A(1)$ on $X$ to algebraic $\tau$-sheaves over $A/\Fa$ on $X$, denoted 
by $\UCF\mapsto \UCF/\Fa\UCF$. If $\UCF$ is a locally free rigid 
analytic $\tau$-sheaf over $A(1)$, then $\UCF/\Fa\UCF$ is a 
locally free algebraic $\tau$-sheaf over~$A/\Fa$. 
 
\bigskip 
 
\noindent{\em Base change\/}: 
 
Let $\UCF$ be a rigid analytic $\tau$-sheaf over $\CA$ on $X$. Further let $\pi:Y\to X$ be a general morphism of rigid analytic spaces $Y$ over $L'$ and $X$ over $L$ (cf.~Definition~\ref{DefGenMorph}). Let 
\[ 
\CA' \enspace := \enspace \CA\widehat{\otimes}_L L' \enspace = \enspace \Gamma(\mathfrak{C}',\CO_{\mathfrak{C}'})\,, \qquad \text{where} \quad \mathfrak{C}' \enspace := \enspace  \mathfrak{C}\widehat{\otimes}_L L'\,. 
\] 
Then we define the rigid analytic $\tau$-sheaf $\pi^\ast\UCF$ over $\CA'$ on $Y$ to be the sheaf $(\pi\times \id_\mathfrak{C})^\ast\CF$ together with the composite homomorphism 
\[ 
\sigma_{Y,\mathfrak{C}'}^\ast \, (\pi\times \id_\mathfrak{C})^\ast\CF \enspace = \enspace (\pi\times \id_\mathfrak{C})^\ast \,\sigma_{X,\mathfrak{C}}^\ast \, \CF \enspace \xrightarrow{\;(\pi\times \id_\mathfrak{C})^\ast\tau_\CF\;} \enspace (\pi\times \id_\mathfrak{C})^\ast\CF\,. 
\] 
In particular, every analytic point $x$ of $X$ (Def.\ \ref{AnalyticPoint}) gives rise to a general 
morphism $i:\Spm k(x) \to X$ (cf.~Section~\ref{AnalyticPoints}). The 
pullback $i^\ast\UCF$ is a $\tau$-sheaf on $\Spm k(x)$, which we also 
denote by $x^\ast\UCF$ in the 
following. For local sections $f$ of $\CF$ we use the notation $f(x) 
:= (i\times \id_\mathfrak{C})^\ast \,f$ for the value of $f$ 
in~$i^\ast\CF$. 
 
\bigskip 
 
{\em $\tau$-invariants I\/}: 
 
A rigid analytic $\tau$-sheaf $\UCF$ over $\CA$ on $X$ comes with a 
homomorphism $\tau_\CF:\enspace\sigma_{X,\mathfrak{C}}^\ast \,\CF 
\longrightarrow \CF$. By adjunction between 
$\sigma_{X,\mathfrak{C}}^\ast$ and $(\sigma_{X,\mathfrak{C}})_\ast$, 
we obtain a homomorphism 
\[ 
 \CF \xrightarrow{\;(\sigma_{X,\mathfrak{C}})_\ast\tau_\CF} (\sigma_{X,\mathfrak{C}})_\ast\,\CF\,, 
\] 
which by abuse of notation we also call $\tau_\CF$. For every {\'e}tale morphism $Y \to X$ it induces a homomorphism on global sections 
\[ 
\tau_\CF: \enspace \Gamma(Y\times_L \mathfrak{C}, \CF) \longrightarrow \Gamma(Y\times_L \mathfrak{C}, (\sigma_{X,\mathfrak{C}})_\ast\,\CF) = \Gamma(Y\times_L \mathfrak{C}, \CF)\,. 
\] 
 
\begin{Def} 
The {\em {\'e}tale sheaf on $X$ of $\tau$-invariants of $\UCF$\/}, 
 denoted by $\UCF^\tau$, is defined as the kernel of 
\[ 
(\id -\tau_\CF):\es W(pr_\ast \,\CF) \longrightarrow W(pr_\ast \,\CF)\,, 
\] 
where $pr:X\times_L \mathfrak{C} \to X$ is the projection onto the first factor. 
I.e., for every {\'e}tale morphism $Y \to X$ we have  
\[ 
\UCF^\tau (Y) \enspace := \enspace \{f \in \Gamma(Y\times_L \mathfrak{C}, \CF): \; \tau_\CF(f) = f\}\,. 
\]  
\end{Def} 
 
The sheaf $\UCF^\tau$ is an {\'e}tale sheaf of $\ul{A}_X$-modules (Sect.~\ref{EtaleSheaves}). If $Y$ is connected and $\UCF$ is 
locally free, we will see in Corollary~\ref{CritUnif} that the sections 
$\UCF^\tau(Y)$ form a projective $A$-module of rank at most~$\rk \CF$. 
To simplify our notation, in later sections we will usually abbreviate 
$W(pr_\ast \CF)$ by~$W(\CF)$. 
 
\bigskip 
 
{\em $\tau$-invariants II\/}: 
 
Suppose $\UCF$ is an algebraic $\tau$-sheaf over $A/\Fa$ on $X$ for 
some non-zero ideal $\Fa$ of $A$. As in the previous construction, one 
can associate to $\UCF$ the {\em {\'e}tale sheaf $(\UCF/\Fa\UCF)^\tau$ of 
$\tau$-invariants modulo $\Fa$}, defined as the kernel of the induced 
homomorphism  
\[ 
(\id -\tau_\CF):\es W(pr_\ast \,\CF) \longrightarrow W(pr_\ast \,\CF)\,. 
\] 
 
The sheaf $\UCF^\tau$ is an {\'e}tale sheaf of $\ul{A/\Fa}_X$-modules.  
 

\section{Torsion Points}\label{TorsSect} 
\setcounter{equation}{0}

Throughout this section let $X$ be a rigid analytic space, $\UCF$ be a 
locally free rigid analytic $\tau$-sheaf over $A(1)$ on $X$ of rank 
$r$, and $\Fa\subsetneq A$ be a non-zero ideal. The $\tau$-sheaf 
$\UCF/\Fa\UCF$ of construction III of change of coefficients  
is locally free on $X\times_{\BF_q} A/\Fa$. By abuse of notation, we 
denote its direct image on $X$ again by $\UCF/\Fa\UCF$. It is locally 
free of rank $r\cdot \dim_{\BF_q} A/\Fa$. On $X$, the homomorphism 
$\tau$ induces a homomorphism 
\begin{equation} \label{CondTauOnF/aF} 
\sigma_{X}^\ast \bigl(\CF/\Fa\CF\bigr) \enspace \longto \enspace \CF/\Fa\CF\,.  
\end{equation} 
 
In later sections we will study the sheaf of global $\tau$-invariant 
sections of $\UCF$, at least {\'e}tale locally. Now  any given section 
of $\UCF^\tau$ will be $\Fp$-adically approximated by a suitable 
compatible systems of sections of the sheaves $(\UCF/\Fp^n\UCF)^\tau$. 
The former element may be thought of as an analytic solution of an 
algebraic equation, the latter as a formal solution to the same equation. The 
formal solutions will have a simpler structure, and in 
general there will be more such than analytic ones. In this section we 
will study formal solutions (and those only up to a finite level), as a 
preparation to the study of analytic solutions in the following 
sections. 
 
\medskip 
 
\begin{Def} 
\begin{enumerate} 
\item 
(cf.\ Taguchi, Wan \cite[{\S} 6]{taguchi-wan} or Drinfeld \cite[{\S} 1]{drinfeld})  
The {\em space of $\Fa$-torsion points of $\UCF$\/} is the rigid analytic group-space on $X$ given by 
\[ 
{}_\Fa \UCF \enspace :=\enspace {\BSpec}_X \, Sym_{\CO_X}(\CF/\Fa\CF)\;/\;(\Frob-\tau) 
\] 
where the ideal $(\Frob-\tau)$ is generated by the local sections $f^q 
-\tau\,f$ for $f\in \CF/\Fa\CF\subset Sym_{\CO_X}(\CF/\Fa\CF)$. The 
structure of commutative group is induced from the group structure of 
$\CF/\Fa\CF$. The action of $A/\Fa$ on $\CF/\Fa\CF$ induces an action $A/\Fa \to \End_X({}_\Fa \UCF)$.  
\item 
The {\em {\'e}tale sheaf $\UCF[\Fa]$ on $X$ of $\Fa$-torsion points of $\UCF$\/} is the sheaf represented by ${}_\Fa \UCF$, i.e. for $\pi:Y\to X$ {\'e}tale we have 
\begin{eqnarray*} 
\UCF[\Fa](Y) & := & \Hom_X(Y,{}_\Fa \UCF) \\[0.2cm] 
& = & \bigl\{ h\in \Hom_{\CO_Y}\bigl(\pi^\ast(\CF/\Fa\CF), \CO_Y\bigr) : h(\tau f)=h(f)^q \\ 
& & \qquad\qquad\qquad\qquad \text{for all local sections } f \text{ of }\pi^\ast(\CF/\Fa\CF)\bigr\} 
\end{eqnarray*} 
The actions of $A/\Fa$ on ${}_\Fa \UCF$ or on $\UCF/\Fa\UCF$ induce 
the same $A/\Fa$-module structures on~$\UCF[\Fa]$. 
\end{enumerate} 
\end{Def} 
 
\begin{remark} 
The group space ${}_\Fa \UCF$ is finite {\'e}tale over $X$ of degree 
$\#(A/\Fa)^r$. This may be verified locally on an admissible affinoid 
covering $\{\Spm B_i\}$ which trivializes $\CF/\Fa\CF$. Let 
$f=(f_1,\ldots,f_s)$ be an $\BF_q$-basis of $\CF/\Fa\CF$ over $\Spm 
B_i$ where $s=r\dim_{\BF_q}A/\Fa$. There is a matrix $\Delta\in 
\GL_s(B_i)$ such that $\tau f=f\Delta$. Then 
\[ 
{}_\Fa \UCF \times_X \Spm B_i \es = \es \Spm B_i[f]/({}^{\sigma\!}f-f\Delta) 
\] 
is finite {\'e}tale over $\Spm B_i$. 
\end{remark}

In addition to these sheaves, we consider the constant sheaf  
$\CH om_{\BF_q}\bigl(\,\nul{A/\Fa}, \BF_q\bigr)$ 
which is fiber-wise free over $A/\Fa$ of rank one. 
 
\begin{Def}\label{DefPerfPairing} 
A pairing $\CG_0\times \CG_1 \to \CH$ of {\'e}tale sheaves $\CG_0,\CG_1,\CH$ of $\nul{A/\Fa}$-modules on $X$ is called {\em perfect} if it induces isomorphisms of $\nul{A/\Fa}$-module sheaves $\CG_i \cong \CH om_{\nul{A/\Fa}}(\CG_{1-i},\CH)$. 
\end{Def}

\begin{lemma}\label{LemmaTorsPoint} 
Suppose the morphism $\tau$ in $(\ref{CondTauOnF/aF})$ is bijective. Then 
\begin{enumerate} 
\item There exists 
a locally finite {\'e}tale covering $\pi\!:Y\to X$ such that  
$(\UCF/\Fa\UCF)^\tau(Y)$ is (on connected components) free over $A/\Fa$ of rank~$r$ and such that 
there is an isomorphism of coherent sheaves 
  $$(\UCF/\Fa\UCF)^\tau(Y)\otimes_{\BF_q}\CO_Y\stackrel{\cong}\longto 
 \pi^*(\CF/\Fa\CF).$$ 
%
\item The pairing  
\begin{equation} \label{ThePairing} 
\UCF[\Fa] \times (\UCF/\Fa\UCF)^\tau \longrightarrow \CH om_{\BF_q\dash\rm Mod}\bigl(\,\nul{A/\Fa}, \BF_q\bigr)\;,\qquad 
(h,f) \longmapsto \bigl( a\mapsto h(af) \bigr)\,. 
\end{equation} 
of {\'e}tale sheaves of ${A/\Fa}$-modules on $X$ is perfect. 
\end{enumerate} 
\end{lemma} 
 
\begin{proof} 
(a) We follow closely Anderson~\cite[Lemma 1.8.2]{anderson}. If $\tau$ 
  is bijective on $\CF/\Fa\CF$, then it is so on 
  $\CF/\Fa^n\CF$. Because some power of $\Fa$ is principal, say equal 
  to $(a)$, it suffices to treat this case. We then choose 
  $\BF_q[t]\subset A$ as the homomorphism defined by $t\mapsto a$. Using 
  Remark~\ref{Taguchi's Trick} we also choose an admissible covering of $X$ by affinoid subdomains such that on each of these subdomains $\Spm B_i$ the sheaf $\CF$ is associated to a free $B_i\langle t\rangle$-module. Then $(\CF/\Fa\CF)(\Spm B_i)$ is a free $B_i$-module of rank $s = r\cdot \dim_{\BF_q}A/\Fa$ on which $\tau$ acts as $\Delta\cdot\sigma$ for a matrix $\Delta\in\GL_s(B_i)$. We consider the equation 
\[ 
\Phi=\Delta {}^{\sigma\!}\Phi \qquad\text{and}\qquad 1=\det \Delta \cdot(\det \Phi)^{q-1} 
\] 
for an $s\times s$-matrix $\Phi$. Let $Y_i:=\Spm B'_i$ be the finite 
{\'e}tale covering of $\Spm B_i$ obtained by adjoining the entries of 
$\Phi$ to $B_i$. Then the columns of the matrix $\Phi\in \GL_s(B'_i)$ 
form a $B'_i$-basis of $(\CF/\Fa\CF)(Y_i)$ as well as a $\BF_q$-basis 
of $(\UCF/\Fa\UCF)^\tau(Y_i)$. So on every connected component of $Y$ 
the $\BF_q$-dimension of $(\UCF/\Fa\UCF)^\tau(Y_i)$ is $r\, 
\dim_{\BF_q}A/\Fa$. Since by \cite[Lemma 1.8.2]{anderson} the stalks 
of $(\UCF/\Fa\UCF)^\tau$ are free $A/\Fa$-modules of 
rank $r$, the proof of (a) is complete. 
 
(b) The question is local on $X$, and so we may assume that the 
isomorphism (a) holds over $X$. Then 
\begin{eqnarray*} 
\UCF[\Fa](X) & = & \bigl\{ h\in \Hom_{\CO_{X}}\bigl(\CF/\Fa\CF, \CO_{X}\bigr) : h(\tau f)=h(f)^q \quad \text{for all }f\bigr\}\\[0.2cm] 
& = & \Hom_{\nul{\BF_q}(X)}\bigl(\,(\UCF/\Fa\UCF)^\tau(X)\,,\,\nul{\BF_q}(X)\,\bigr)\\[0.2cm] 
& \cong & \Hom_{\nul{A/\Fa}(X)}\bigl(\,(\UCF/\Fa\UCF)^\tau(X)\,,\,  
\CH om_{\BF_q}\bigl(\,\nul{A/\Fa}, \BF_q\bigr)(X)\,\bigr)\,, 
\end{eqnarray*} 
the last isomorphism mapping $h\in \UCF[\Fa](X)$ to $f\longmapsto 
\bigl(a\mapsto h(af)\bigr)$. The same chain of isomorphisms holds for 
all {\'e}tale morphisms $Y\to X$, and so the pairing is perfect. 
\end{proof}

\begin{thm} \label{ThmaTorsion} 
Suppose that $(\ref{CondTauOnF/aF})$ is an isomorphism. Then there exists a (finite) {\'e}tale $\GL_r(A/\Fa)$-torsor $\pi:X_{\UCF,\Fa}\to X$ such that 
\begin{enumerate} 
\item 
$\UCF[\Fa](X_{\UCF,\Fa})$ is a free $\nul{A/\Fa}(X_{\UCF,\Fa})$-module of rank $r$ and 
\item 
$(\UCF/\Fa\UCF)^\tau(X_{\UCF,\Fa})$ is a free $\nul{A/\Fa}(X_{\UCF,\Fa})$-module of rank $r$. Furthermore, the natural homomorphism of sheaves on $Y$ 
\[ 
(\UCF/\Fa\UCF)^\tau(X_{\UCF,\Fa}) \otimes_{\BF_q} \CO_{X_{\UCF,\Fa}} \rbij \pi^\ast\,\CF/\Fa\CF 
\] 
is an isomorphism. 
\end{enumerate} 
\end{thm}

\begin{proof} 
(a) This follows from general principles on sections of finite {\'e}tale group schemes; cf.\ Katz-Mazur~\cite[Chap.\ 1]{Katz-Mazur}. To be precise ${}_\Fa \UCF$ is a finite {\'e}tale commutative group over $X$ of rank $(\# A/\Fa)^r$. Thus the functor which assigns to every rigid analytic $X$-space $T$ the set 
\[ 
\nul{\Isom}_{X\text{-Gp}}(\bigl(\nul{A/\Fa}\bigl)^{ r}, {}_\Fa \UCF)\,(T) \enspace = \enspace \{\text{group isomorphisms}\quad \bigl(\nul{A/\Fa}\bigl)^{ r}(T) \rbij {}_\Fa \UCF(T)\,\} 
\] 
is representable by a rigid analytic space $Y'$, which is finite {\'e}tale over $X$; cf.\ \cite[Cor.\ X.5.10]{SGA3}. 
(To apply the theory written for schemes to rigid analytic spaces, cover $X$ by affinoids $\Spm B_i$ and work over $\Spec B_i$. The representing schemes $Y_i$ are finite {\'e}tale over $\Spec B_i$ and induce rigid analytic spaces which glue due to their universal property.) Therefore the open and closed sub-functor  
\begin{eqnarray*} 
&&\nul{\Isom}_{\nul{A/\Fa}\text{-Mod}}(\bigl(\nul{A/\Fa}\bigl)^{ r}, {}_\Fa \UCF)\,(T) \enspace = \\ 
&& \qquad \qquad = \enspace \{A/\Fa-\text{module isomorphisms}\quad \bigl(\nul{A/\Fa}\bigl)^{ r}(T) \rbij {}_\Fa \UCF(T)\,\}\enspace \\ 
&& \qquad \qquad \subset \enspace \nul{\Isom}_{X\text{-Gp}}(\bigl(\nul{A/\Fa}\bigl)^{ r}, {}_\Fa \UCF)\,(T) 
\end{eqnarray*} 
is representable by a rigid analytic space $Y\subset Y'$, which is finite {\'e}tale over $X$. Above every connected component of $X$ the space $Y$ is either empty or a $\GL_r(A/\Fa)$-torsor. Finally the surjectivity of $\pi:Y\to X$ follows from Lemma~\ref{LemmaTorsPoint}. 
 
(b) We deduce from (a) and Lemma~\ref{LemmaTorsPoint} that 
\[ 
(\UCF/\Fa\UCF)^\tau(X_{\UCF,\Fa}) \es \cong \es \CH om_{\nul{A/\Fa}}\bigl(\UCF[\Fa], \CH om_{\BF_q}(\nul{A/\Fa},\BF_q)\bigr)\,(X_{\UCF,\Fa}) 
\] 
is a free $\nul{A/\Fa}(X_{\UCF,\Fa})$-module of rank $r$. The last statement also follows from Lemma~\ref{LemmaTorsPoint}. 
\end{proof}

 
\section{The Key-Lemma for $\tau$-Sheaves} 
\setcounter{equation}{0}\label{KeyLemSec}

In this section we prove the following Key-Lemma on the extension of sections of a $\tau$-sheaf. We will apply the Key-Lemma to show that $\UCF^\tau$ is an overconvergent sheaf (Def.~\ref{DefOverconvergent}). In the following sections the Key-Lemma will serve as the main tool for studying the analytic triviality of $\tau$-sheaves. For the required background on analytic points see the appendix \ref{AnalyticPoints}.

\begin{keylemma} \label{ExtendOnV} 
Let $X = \Spm B$ be affinoid, set $\FC:=\mathfrak{A}(1)$ and let $\UCF$ be a rigid analytic $\tau$-sheaf over $L\langle t \rangle$ on $X$ such that $\CF$ is a free $B \langle t \rangle$-module of rank $r$. Let $g_1, \ldots, g_m \in \Gamma\bigl(X \times_L \mathfrak{C},\CF\bigr)$. Let further $x \in \CM(X)$ be an analytic point and let $\ol{f}_1, \ldots, \ol{f}_m$ be given elements of $x^\ast\CF)\bigl(k(x)^{alg}\bigr)$, defined over an algebraic closure of $k(x)$, which satisfy $\ol{f}_i - \tau\,\ol{f}_i = g_i(x)$ for all $i$. Then: 
 
\begin{enumerate} 
\item 
There exists an {\'e}tale morphism $\pi:V \to X$ of affinoids such that $\pi V$ is a wide neighborhood of $x$ in $X$, a point $y\in \CM(V)$ above $x$, and uniquely determined sections $f_1, \ldots, f_m \in \Gamma\bigl(V \times_L \mathfrak{C},\CF\bigr)$ satisfying $f_i(y) = \ol{f}_i$ and $f_i - \tau\,f_i = \pi^\ast g_i$ for all $i$. Moreover, 
\item If $\tau$ is bijective on $\CF/t\CF$, then $V$ may be taken as an 
affinoid subdomain of $X_{\UCF,t^N}$, constructed in 
Theorem~\ref{ThmaTorsion} (with the obvious morphism to $X$) for some 
$N\gg0$. 
\item 
If all $\ol{f}_i$ are defined over $k(x)$, we can find $V \subset X$ as a wide affinoid neighborhood of $x$. 
\end{enumerate} 
\end{keylemma} 
 
\medskip 
 
\begin{proof} 
We have to find $m$ sections of $\CF$, which we will represent by an $r \times m$ matrix $\Phi$, using a basis of $\CF$ over $B\langle t \rangle$. 
We choose this basis such that $\sup_n|\Delta_n|= 1$, cf.\ Remark~\ref{Taguchi's Trick}. Let $\Omega = \sum_{n\geq 0} \Omega_n t^n \in M_{r\times m}(B\langle t \rangle)$ be the coordinate matrix of the system $(g_1, \ldots, g_m)$ with respect to this basis. 
The condition $f_i - \tau\,f_i = g_i$ translates into the system of equations 
\begin{equation} \label{EqTauInv} 
\Phi_n - \,\Delta_0 \,{}^{\sigma\!} \Phi_n \enspace = \enspace \Psi_n \enspace := \enspace \Omega_n \,+\,\sum_{\nu =1}^n \Delta_\nu \,{}^{\sigma\!} \Phi_{n-\nu}\,, \qquad n \geq 0\,. 
\end{equation} 
 
For each of the matrices $\Phi_n$ this is an {\'e}tale equation. Note that if $|\Psi_n|<1$, then there is a unique solution $\Phi_n$ with $|\Phi_n|= |\Psi_n|$, given explicitly by 
\begin{equation} \label{ASSolution} 
\Phi_n \enspace = \enspace  \sum_{\nu=0}^\infty \,\Delta_0 \;{}^{\sigma\!}\Delta_0 \cdot \ldots \cdot{}^{\sigma^{\nu-1}}\!\Delta_0\;{}^{\sigma^\nu}\!\Psi_n\,. 
\end{equation} 
 
Set $\DS(\ol{f}_1, \ldots, \ol{f}_m) = \sum_{n\in \Nz_0} \ol{\Phi}_n t^n =: \ol{\Phi}$ with matrices $\ol{\Phi}_n \in M_{r \times m}\bigl(k(x)\bigr)$.  
 
If $x$ is a classical point and the $\ol{f}_1, \ldots, \ol{f}_m$ are defined over $k(x)$, then we can use the Inverse Function Theorem to extend $\ol{\Phi}$ to a wide neighborhood $U$ of $x$. Namely, we lift $\ol{\Phi}$ to an element of $M_{r\times m}(B\langle t\rangle)$ and set $\Phi = \ol{\Phi}+ \Phi'$. If $U$ is such that $|\Phi'_\nu|\ll 1$ on $U$ for all $0\leq \nu <n$, there is a unique solution $\Phi'_n$ on $U$ with $|\Phi'_n|\ll 1$ to equation~(\ref{EqTauInv}) as in (\ref{ASSolution}).  
 
In the general case however, we can not explicitly solve equation~(\ref{EqTauInv}) for all $n$. Instead, we formally adjoin the solutions of the equations for $n=0, \ldots, N$, thus obtaining an {\'e}tale morphism $V \to X$. On $V$ we now may solve equation~(\ref{EqTauInv}) explicitly for all $n>N$ as in (\ref{ASSolution}).  
 
To be precise, we fix constants $\Theta$ and $\epsilon$ in $|L^{\alg}|$ satisfying 
\begin{eqnarray*} 
0 \enspace < & \Theta & <  \enspace {\TS \frac{1}{2}} \qquad\text{and}\\[0.2cm] 
0 \enspace < & \epsilon & < \enspace \inf \,\{\;1\,,\enspace {\TS \frac{1}{2}}\,\Theta \,|\ol{\Phi}_n|^{-q}: n \geq 0\,\} \enspace \leq \enspace 1\,. 
\end{eqnarray*} 
 
There is an $l\in \Nz$ such that  
\[ 
|\Delta_n| \enspace \leq \enspace \epsilon \quad, \qquad |\ol{\Phi}_n| \enspace \leq \enspace {\TS \frac{1}{2}}\,\Theta \qquad \text{and} \qquad |\Omega_n| \enspace \leq \enspace {\TS \frac{1}{2}}\,\Theta \qquad \text{for every} \quad n>l\,. 
\] 
 
We set $N :=2l$ and define $B_N$ to be the polynomial ring over $B$ in the components of the matrices $\Phi_0, \ldots, \Phi_N$ modulo the relations (\ref{EqTauInv}) for $n=0, \ldots, N$.  
 
If $\tau$ is bijective on $\CF/t\CF$, the matrix $\Delta_0$ is 
invertible. In this case $B_N$ is an affinoid algebra and $\pi:Y:=\Spm 
B_N \to \Spm B$ is the finite {\'e}tale $\GL_r(\BF_q[t]/t^N)$-torsor above $\Spm B$ which trivializes the $t^{N+1}$-torsion of $\UCF$; cf. Thm.\ \ref{ThmaTorsion}.  
Now consider the affinoid subsets 
\begin{eqnarray*} 
V & := & \{y \in \Spm B_N: |\Phi_n(y)| \leq \Theta \quad \forall \,l < n \leq N \quad \text{and} \quad |\Phi_n(y)|^q \leq \Theta/\epsilon \quad \forall \,0\leq n \leq l \}\\[0.2cm] 
U & := & \{y \in \Spm B_N: |\Phi_n(y)| \leq \frac{\Theta}{2}\quad \forall \,l < n \leq N \quad \text{and} \quad |\Phi_n(y)|^q \leq \frac{\Theta}{2\epsilon} \quad \forall \,0\leq n \leq l \} 
\end{eqnarray*} 
of $\Spec B_N$. 
 
The choice of our constants implies that $x$ lifts to an analytic point $y \in \CM(U)$ with $\ol{\Phi}_n = \Phi_n(y)$ for all $n\leq N$. We claim that equation (\ref{EqTauInv}) possesses solutions $\Phi_n$ in $\Gamma(V,\CO_V)$ for all $n >N$ which satisfy $|\Phi_n|\leq \Theta$. Indeed, by induction  
\[ 
|\Psi_n| \enspace \leq \enspace \max \{ \,|\Omega_n|\,,\enspace|\Delta_\nu|\, |\Phi_{n-\nu}|^q\;,\;1\leq \nu \leq n\,\} 
\] 
and for $n>N$ the right hand side is at most the maximum of $|\Omega_n| \leq \frac{1}{2}\Theta$ and 
\begin{eqnarray*} 
& & \sup \{\,|\Delta_n|: n \geq 0\,\} \cdot \Theta^q \enspace = \enspace \Theta^q\enspace \leq \enspace {\TS \frac{1}{2}}\,\Theta \qquad \text{for} \quad \nu \leq l\\[0.2cm] 
& & \epsilon  \cdot \bigl(\sup \{\,|\Phi_n|: n \geq 0\,\}\bigr)^q \enspace \leq \enspace \epsilon \; \Theta/\epsilon \enspace = \enspace \Theta  \qquad \text{for} \quad \nu > l\,. 
\end{eqnarray*} 
 
Hence we may take for $\Phi_n$ the solution to equation~(\ref{EqTauInv}) given by (\ref{ASSolution}). It is the unique solution with $|\Phi_n|= |\Psi_n|\leq \Theta$. Therefore we have $\ol{\Phi}_n = \Phi_n(y)$ also for all $n>N$. The same argument together with $\DS\lim_{\nu \to \infty}|\Delta_\nu|= 0$ and $\DS\lim_{n \to \infty}|\Omega_n|=0$ shows that 
\[ 
\limsup_{n\to \infty} |\Phi_n| \enspace \leq \enspace \Bigl(\limsup_{n\to \infty} |\Phi_n|\Bigr)^q \enspace \leq \Theta^q\enspace<\enspace1, 
\] 
hence $\DS\lim_{n\to \infty}|\Phi_n|=0$. 
So the columns of  
\[ 
\Phi \enspace := \enspace \sum_{n=0}^\infty\, \Phi_n t^n 
\] 
are the desired elements $f_1,\ldots,f_m$ satisfying $f_i - \tau\,f_i = \pi^\ast g_i$. These are uniquely determined by the condition $f_i(y) = \ol{f}_i$, since the defining equation (\ref{EqTauInv}) for the $\Phi_n$ is {\'e}tale.

Let $\wt{\Theta}, \wt{\epsilon} \in L^{\alg}$ be constants with $|\wt{\Theta}|=\Theta$ and $|\wt{\epsilon}|=\epsilon$. Then the elements $\wt{\Theta}^{-1}\Phi_n \in B_N$ for $l < n \leq N$ and $(\wt{\epsilon}/\wt{\Theta})^{1/q} \,\Phi_n \in B_N$ for $0\leq n\leq l$ form an affinoid generating system of $V$ over $X$. Therefore we see that $U$ is relatively compact in $V$ over $X$, i.e., $U \subset \subset_X V$. Thus by~
\cite[Lemma 3.4.2]{JongPut} the image $\pi(V)$ is a wide neighborhood of $\pi(U)$ in $X$. Since $x \in \pi\CM(U)=\CM(\pi U)$ by \cite[Prop 3.1.7]{JongPut}, we conclude that $\pi V$ is a wide neighborhood of $x$ in $X$. 
 
If all $\ol{f}_i$ are defined over $k(x)$, we have $k(y) =k(x)$. By \cite[Lemma 3.1.5]{JongPut} there is an affinoid subdomain $U'\subset U$ with $y\in\CM(U')$ such that $\pi|_{U'}$ is an isomorphism $U'\to \pi(U')$. Now by \cite[Lemma 3.4.2]{JongPut} there is a wide neighborhood $V'$ of $y$ such that $\pi|_{V'}:V'\to \pi V'$ is an isomorphism and $\pi V'$ is a wide neighborhood of $x$ in $X$, i.e., $\CM(\pi V')$ is a neighborhood of $x$ in $\CM(X)$. 
\end{proof}

\medskip 
 
\begin{cor} \label{SatUModTau} 
Let $\UCF$ be a locally free rigid analytic $\tau$-sheaf over $A(1)$ on a connected $X$ and let $\pi :Y\to X$ be a general morphism, cf.~Definition~\ref{DefGenMorph}. Then the homomorphism of $A$-modules 
\[ 
\UCF^\tau(X) \longrightarrow (\pi^\ast\UCF)^\tau(Y)\,,\quad f \mapsto \pi^\ast f 
\] 
is injective and its image is a saturated $A$-submodule. 
\end{cor}

\begin{proof} 
For the proof, we may assume that $A=\BF_q[t]$. Since $\UCF^\tau$ is also a Zariski sheaf,  
we may assume that $X$ is an affinoid over which we have local coordinates for $\UCF$ as  
in Remark~\ref{Taguchi's Trick}. 
 
The proposition is clear in the case when $X$ and $Y$ are spectra of 
algebraically closed complete fields. Let $y$ be an analytic point of $Y$ and 
$x= \pi(y)$. We let $\wh{k(x)^{alg}}$ denote the completion of an algebraic closure of $k(x)$. By considering the diagram 
\[\xymatrix @C+2pc{ 
\UCF^\tau(X) \ar[r]\ar[d] &(\pi^*\UCF)^\tau(Y)\ar[d]\\ 
(x^\ast\UCF)^\tau\bigl(\wh{k(x)^{alg}}\bigr) \ar[r]^\cong & 
(x^\ast\UCF)^\tau\bigl(\wh{k(y)^{alg}}\bigr) 
} 
\] 
the assertion on injectivity is clear, if we have shown it for the vertical 
morphisms, i.e. it remains to consider the case $Y=\Spm \wh{k(x)^{alg}}$. 
In local coordinates the defining equation (\ref{EqTauInv}) for 
the coefficients of $\tau$-invariants are {\'e}tale. Therefore if we have  
a global section which vanishes at $x$ ({\'e}tale locally), it will have to 
vanish on the entire connected component of $X$ containing~$x$. Since $X$ is connected, 
injectivity is established.

To prove saturation of $h$, again by the above diagram 
it suffices to consider the case where $Y=\Spm \wh{k(x)^{alg}}$. Let now  
$f$ be in $\UCF^\tau(X)$, $0\neq \ol{g}$ in  
$(x^\ast\UCF)^\tau\bigl(\wh{k(x)^{alg}}\bigr)$,  
and $a \in A$ such that $f(x) = a\,\ol{g}$. We have to show 
that $\ol{g}$ extends to a section of $\UCF^\tau$ over $X$. We may 
also assume that $a$ is not a constant, since then the assertion is clear.

By replacing $A$ by $\BF_q[a]$, we may assume $a=t$. 
Let $\sum_n v_nt^n$ be the series representing 
$\bar g$ with coefficients in $\wh{k(x)^{alg}}$. Then $\sum v_nt^{n+1}$ 
represents $f$ and has coefficients in $\Gamma(X,\CO_X)$. But then the 
$v_n$ are in $\Gamma(X,\CO_X)$, and hence $\sum_nv_nt^n$ lies in 
$\UCF^\tau(X)$ and maps to~$\bar g$.  
\end{proof}

\begin{cor}\label{StalkOfFtauAndOverC} 
Let $\UCF$ be a locally free rigid analytic $\tau$-sheaf over 
$A(1)$. Then 
\begin{enumerate} 
\item \label{StalkOfFtauAndOverC_Parta} 
For every analytic point $x\in\CM(X)$ the {\'e}tale sheaves $x^\ast(\UCF^\tau)$ and $(x^\ast\UCF)^\tau$ on $\Spm k(x)$ are canonically isomorphic. 
\item\label{StalkOfFtauAndOverC_Partb} 
For every {\'e}tale analytic point $y\in \CM_\et(X)$ (Sect.~\ref{EtaleSheaves}), the stalk of $\UCF^\tau$ at $y$ is given by 
  $(y^\ast\UCF)^\tau\bigl(k(y)\bigr)$. 
\item \label{StalkOfFtauAndOverC_Partc} 
$\UCF^\tau$ is an overconvergent sheaf on the small {\'e}tale site 
  of~$X$. 
\item \label{StalkOfFtauAndOverC_Partd} 
The homomorphism $(\id-\tau):W(\CF)\to W(\CF)$ is surjective 
  if and only if it is so for the fibers $y^\ast\UCF$ at all {\'e}tale analytic  
  points $y\in \CM_\et(X)$. 
\end{enumerate} 
\end{cor}

\begin{proof} 
All assertions are local. Thus by Lemma~\ref{lubolemma} and 
Remark~\ref{Taguchi's Trick} we may assume that $A=\BF_q[t]$ and that 
$\CF$ is free of finite rank on $X\times \FA(1)$.  To see part~\ref{StalkOfFtauAndOverC_Parta}, note that by  
the previous corollary the canonical morphism of {\'e}tale sheaves $x^\ast(\UCF^\tau)\to (x^\ast\UCF)^\tau$ is injective. That it surjective follows directly from  
Lemma~\ref{ExtendOnV}~(a). Next, part~\ref{StalkOfFtauAndOverC_Partb} is just a reformulation of \ref{StalkOfFtauAndOverC_Parta}. 
 
For \ref{StalkOfFtauAndOverC_Partc} let $Y\to X$ be {\'e}tale and $V\subset Y$ a special 
subset. We need to show that 
$$\UCF^\tau(V)=\mathop{\dirlim}\limits_{V\subset\subset_YU}\UCF^\tau(U),$$  
and so let $f \in \UCF^\tau(V)$.  By Lemma~\ref{ExtendOnV}~(c) with $m=1$ 
and $g_1=0$, there exists for every analytic point $y$ of $V$ a wide 
affinoid neighborhood $U_y$ of $y$ in $Y$ and a unique section 
$\widetilde{f}_y \in \UCF^\tau(U_y)$ such that $f(y)=\wt f_y(y)$. By the 
compactness of $\CM(V)$ finitely many of the $\CM(U_y)$ will cover 
$\CM(V)$. So there are finitely many $y$ such that the special subset 
$U := \cup \,U_y$ is a wide neighborhood of $V$ in $Y$. Since the 
sections $\widetilde{f}_y$ agree with $f$ on $V$, by the previous  
corollary applied to suitable points in $V$ and to $Y$, they glue  
to a section $\widetilde{f} \in \UCF^\tau(U)$. 
 
Finally \ref{StalkOfFtauAndOverC_Partd} is a direct consequence of the first assertion of the Key-Lemma~\ref{ExtendOnV}~(a), since if a special set $X$ is covered by wide 
open sets, it is covered by finitely many of these. 
\end{proof}

\smallskip

\begin{prop} \label{TauInvOnAInfty} 
Let $\UCG$ be a locally free $\tau$-sheaf over $A(\infty)$ 
on $X$, and let $\UCF$ be the locally free $\tau$-sheaf over $A(1)$ on 
$X$, obtained from $\UCG$ by the change of coefficients II (Sect.\ \ref{Functors}). Then the natural 
homomorphism $\UCG^\tau \to \UCF^\tau$ is an 
isomorphism of {\'e}tale sheaves of $A$-modules. 
\end{prop} 
 
\begin{proof} 
We may assume that $A=\BF_q[t]$. Choose $c\in L$ such that $|c|>1$ and 
define $D_m:=\Spm L\langle \frac{t}{c^m}\rangle$, so that $X \times_L 
D_m$, $m\in\BN$, is an admissible covering of $X \times_L (\Spec 
A)^{\rig}$. We have to show that $(\UCG|_{X\times_L D_m})^\tau\to\UCF^\tau$ is 
an isomorphism for all~$m\geq0$.  
 
Since both sheaves are overconvergent 
by Corollary~\ref{StalkOfFtauAndOverC}, it suffices to verify the 
isomorphy on 
{\'e}tale stalks, and hence to assume that $X=\Spm F$ for some 
algebraically closed complete normed field $F$ above $L$. The 
$\tau$-sheaf $\UCG$ is thus given by a free finitely generated module $M$ over  
$$F\langle \langle t\rangle \rangle:=\Big\{\sum_{\nu\geq0} a_\nu t^\nu\bigm| 
a_\nu\in F;\,\forall\rho>1:\lim_{\nu\to\infty}|a_\nu|\rho^\nu=0 \Big\}$$  
of some rank $r$, on which $\tau_\CG$ is represented by a matrix 
$\Delta = \sum \Delta_\nu t^\nu \in M_r\bigl(F\langle\langle t\rangle 
\rangle\bigr)$, i.e., 
\begin{equation}\label{LimForDel} 
\forall\rho>1:\lim_{\nu \to \infty} |\Delta_\nu|\,\rho^\nu \enspace = \enspace 0\,. 
\end{equation} 
The elements of $\UCF^\tau(X)$ can be written as vectors $\Phi=\sum 
\Phi_n t^n \in F\langle t\rangle^r$, i.e., with 
\begin{equation}\label{LimForPhi} 
\DS\lim_{n\to\infty}|\Phi_n| \es=\es 0 
\end{equation} 
and subject to the condition $\Phi = \Delta\,{}^{\sigma\!}\Phi$, i.e., 
\[ 
\Phi_n - \,\Delta_0 \,{}^{\sigma\!} \Phi_n \enspace = \enspace \Psi_n 
:= \enspace \sum_{\nu =1}^n \Delta_\nu \,{}^{\sigma\!} \Phi_{n-\nu}\,, \quad n \geq 0\,. 
\] 
 
We set $\rho:=c^{m+1}$, $m\in\BN$ and claim that $|\Phi_n|\rho^n$ 
is bounded -- note that this implies that $|\Phi_n||c^m|^n$ converges to zero. As $m$ was arbitrary this will finish the proof. 
Now the sequences in (\ref{LimForDel}) and (\ref{LimForPhi}) are zero sequences. Therefore 
we may find some $N\in\BN$ such that for all $n > N$ and 
$\nu\in\{0,\ldots,n\}$ we have 
\[ 
|\Delta_\nu| \,\rho^\nu\,|\Phi_{n-\nu}|^{q-1} \enspace \leq \enspace 1/2 \,. 
\] 
Let now $n>N$. For $\nu=0$ the inequality implies $|\Delta_0 
\,{}^{\sigma\!}\Phi_n|<|\Phi_n|$, and so $|\Phi_n|=|\Psi_n|$. Thus 
\begin{eqnarray*} 
|\Phi_n|\,\rho^n & \leq &  \rho^n\max_{\nu=1,\ldots,n}  |\Delta_\nu| |\Phi_{n-\nu}|^q 
\enspace \leq \enspace  
\max_{\nu=1,\ldots,n}  \Big(|\Delta_\nu| \rho^\nu |\Phi_{n-\nu}|^{q-1}\Big) |\Phi_{n-\nu}|\rho^{n-\nu}\\ 
& \leq & 1/2 \max_{\nu=0,\ldots,n-1} |\Phi_{\nu}|\rho^{\nu}, 
\end{eqnarray*} 
which shows that $|\Phi_n|\rho^n$ is bounded and thus proves the claim. 
\end{proof}


\section{Triviality of $\tau$-Sheaves} \label{SectTriv} 
\setcounter{equation}{0}

The simplest $\tau$-sheaf is $\BOne_{X,A(1)}$. This section is 
centered around the question of how to recognize whether a given 
(locally free) rigid analytic $\tau$-sheaf is (basically) of this 
trivial form. In the following section, we study the same question {\'e}tale locally 
which leads to the notion of analytic triviality.  
 
\begin{Def}\label{TrivialDef} 
A (locally free) rigid analytic $\tau$-sheaf $\UCF$ on $X$ is \emph{trivial} 
if there exists a finitely generated projective $A$-module $P$ such 
that $\UCF\cong P\otimes_A\BOne_{X,A(1)}$. 
\end{Def} 
If $\UCF$ is trivial, then clearly $P\cong \UCF^\tau(X)$. 
 
\medskip 
 
Let us first rephrase some results of Anderson, \cite[{\S}~2]{anderson}, 
in our language. 
\begin{lemma}\label{AndersRephr} 
Suppose $L'$ is a finite field extension of $L$ and $X=\Spm L'$.  
\begin{enumerate} 
\item Any sub-$\tau$-sheaf of a trivial $\tau$-sheaf on $X$ is trivial. 
\item Suppose $\UCF$ is a locally free rigid analytic $\tau$-sheaf 
  over $A(1)$ on $X$ with module of global $\tau$-invariants 
  $P:=\UCF^\tau(X)$. Then $P$ is finitely generated projective over 
  $A$ and the induced homomorphism 
  $$h\!: P\otimes_A\BOne_{X,A(1)}\to \UCF$$ 
  is injective. In particular the rank of $P$ over $A$ is bounded 
  by the rank of $\UCF$ over $X\times\FA(1)$. 
\item If either the image of $h$ in (b) is of full rank, or 
  equivalently (by (b)), if $\UCF$ and $P$ are of the same rank, then 
  $h$ is an isomorphism. 
\end{enumerate} 
\end{lemma} 
\forget{Suppose $\UCG\subset\UCF$ is a $\tau$-subsheaf and $\CF$ is locally free. Then by the saturation of $\CG$ in $\CF$ we 
mean the subsheaf $\CF'$ whose module of global sections is the 
intersection of the module of global sections of $\CF$ with the 
generic stalk of~$\CG$. One can easily show that $\tau$ preserves 
$\CF'$, and so the saturation $\UCF'$ is again a $\tau$-subsheaf.} 
\begin{proof} 
To prove (a), suppose $\UCF=P\otimes_A\BOne_{X,A(1)}$ is trivial, 
suppose $\UCG\subset\UCF$ is a $\tau$-subsheaf  and set 
$P:=\UCF^\tau(X)$ and $Q:=\UCG^\tau(X)$. It will suffice to show that 
$Q\otimes_A\BOne_{X,A(1)}\to \UCG$ is an isomorphism. To prove this, 
we may assume that~$A=\BF_q[t]$.  
Let $\UCF''$ be the cokernel of $\UCG\into\UCF$. 
By Anderson \cite[Lemma 2.10.4]{anderson} there is a finite 
field extension $L''$ of $L'$ such that  
\[ 
\UCG\otimes_{L'} L''\,\cong \,(\UCG\otimes_{L'} L'')^\tau(L'') \,\otimes_A\BOne_{\Spm L'',A(t)}\,. 
\] 
Because the functor of $\tau$-invariants is left 
exact and because of Corollary~\ref{SatUModTau}, we have the 
commutative diagram  
$$ 
\xymatrix @!C=4pc {0\ar[r]&\ar@{^{ (}->}[d]Q\ar[r]&P\ar[r] \ar@{=}[d]& 
  \ar@{^{ (}->}[d] (\UCF'')^\tau(X)\\ 
0\ar[r]&(\UCG\otimes_{L'} L'')^\tau(L'')\ar[r]&P\ar[r] & (\UCF''\otimes_{L'} L'')^\tau(L'')\rlap{,}\\}$$ 
in which the rows are left exact and the vertical homomorphisms are 
injective. Since $\UCF''$ is not necessarily locally free, strictly speaking Corollary~\ref{SatUModTau} does not apply. 
But we may apply it to the right hand term in $0\to \UCF''_\tors\to\UCF''\to\UCF''/\UCF''_\tors$, and on the 
left the injectivity assertion is clear from Section~\ref{TorsSect}.  
 
A simple diagram chase implies that $Q\into (\UCG\otimes_{L'} 
L'')^\tau(L'')$ must be an isomorphism, and so in particular the rank of 
$Q$ over $A$ agrees with the rank of $\CG$ over 
$X\times\FA(1)$. Therefore $Q\otimes_A \BOne_{X,A(1)}\into \UCG$ 
is a $\tau$-subsheaf of full rank. Since after base change along the 
faithfully flat morphism $L'\to L''$ it becomes an isomorphism, it 
must have been an isomorphism to begin with. This proves~(a). 
 
\forget{From the explicit form of the homomorphism $\UCG\otimes_{L'} 
L''\into \UCF\otimes_{L'} L''$, it follows that  
$\UCF''\otimes_{L'} L''\cong  P/P'\otimes_A\BOne_{\Spm L',A(1)}$. 
So considering dimensions in the above diagram, after tensoring with 
$K$ over $A$, we find that the ranks of $Q$ and $P'$ over $A$ must be 
the same. But the again by Corollary~\ref{SatUModTau} we must have 
$Q=P'$. Because $\UCF''$ is locally free, $P/P'$ must be so, too, and 
thus $Q\otimes_A\BOne_{X,A(1)}\subset \CG$ is saturated and of the 
same rank as $\CG$. Therefore the two agree, and this proves (a) for 
saturated~$\UCG$. 
 
In the general case, we may, by the case already treated, assume that 
$\UCF$ is the saturation of $\UCG$. Let $\UCF''=\UCF/\UCG$ as above. 
In the present case $\UCF''$ is of finite length, and hence 
$(\UCF'')^\tau$ is finite. Thus $P$ and $Q$ must be of the same rank. 
Let again $L''/L'$ be a finite extension which trivializes $\UCG$. 
As above we must have $Q\cong (\UCG\otimes_{L'}L'')^\tau(X)$ by 
Corollary~\ref{SatUModTau}. This shows that 
$Q\otimes_A\BOne_{X,A(1)}\to \UCG$ becomes an isomorphism under the 
faithfully flat base change $L'\to L''$, and hence is an isomorphism 
itself.} 
 
\smallskip 
 
To prove (b), we may replace $\UCF$ by the image of $h$, since the 
latter is a locally free $\tau$-subsheaf of $\UCF$ and since it has 
the same $\tau$-invariants on $X$ as $\UCF$. Considering the ranks of the 
domain and range of $h$, there must exist a projective $A$-submodule $P'$ of 
$P$ of rank equal to the rank of $\UCF$. Let $\UCG$ denote the kernel 
of $\UCF':=P'\otimes_A\BOne_{X,A(1)}\to \UCF$. Since taking $\tau$-invariants 
is left exact, by the definition of $h$ we must have 
$\UCG^\tau(X)=0$. By part (a) it follows that $\UCG=0$, and hence $h$ 
is injective on $\UCF'$. Therefore the cokernel of  
$h\!:\UCF'\to\UCF$ must be of finite length over 
$X\times_L\FA(1)$. Taking $\tau$-invariants yields 
$$0\longto P'\longto P\longto (\UCF'')^\tau.$$ 
As before, the term on the right is finite, and this completes the 
proof of~(b). 
 
\smallskip 
 
For part~(c), it remains to prove the surjectivity of $h$. We may 
again assume $A=\BF_q[t]$. Then $P$ is free, and so we may 
choose an $A$-basis $f_1, \ldots, f_r$ of it. We now use that 
$k(x)\langle t\rangle$ is a principal ideal domain, as is shown for 
instance in \cite[{\S}~2]{anderson}, and denote by  $\alpha_\nu$ the 
elementary divisors of a matrix representing $h$. If we apply the 
$r$-th exterior power operation to $h$, the resulting elementary 
divisor is $\alpha_1\cdot\ldots\cdot \alpha_r$. Because $\tau$ 
preserves the image, the product ${}^{\sigma\!}\alpha_1 \cdot 
\ldots\cdot{}^{\sigma\!}\alpha_r$ divides $\alpha_1 \cdot 
\ldots\cdot\alpha_r$, and thus all $\alpha_\nu$ lie in $\Fz_q[t] \cdot 
k(x)\langle t\rangle\mal$. It follows that there is an $a\in 
\Fz_q[t]$ such that $a \CF$ is a $\tau$-submodule of $\im(h)\cong  
P \otimes_A L'\langle t\rangle$. By part (a) we must have 
$a\UCF\cong (a\UCF)^\tau(X)\otimes_A\BOne_{X,A(1)}$, and (c) is now 
immediate. 
\end{proof} 
 
\begin{cor}\label{CritUnif} 
Suppose $X$ is connected and $\UCF$ is a locally free rigid analytic 
$\tau$-sheaf over $A(1)$ on $X$. Then $P:=\UCF^\tau(X)$ is a finitely 
 generated projective $A$-module of rank at most the rank of~$\CF$. 
Equality of ranks holds if and only if $\UCF$ is trivial. 
\end{cor} 
 
\begin{proof} 
By Corollary~\ref{SatUModTau}, the first assertion may be 
reduced to the case where $X$ the spectrum of $k(x)$ with $x$ a classical 
point on~$X$. However this case is proved in the previous lemma.  
It remains to show that 
\[ 
h\!:\BOne_{X,A(1)} \,\otimes_A \UCF^\tau(X) \longrightarrow \UCF\,, \quad b \otimes f \mapsto b\,f\,. 
\] 
is an isomorphism of rigid analytic $\tau$-sheaves on $X \times_L 
\mathfrak{A}(1)$ provided that the ranks of $\CF$ and $P$ agree, since 
the converse implication is obvious. 
 
Because domain and range of $h$ are locally free of the same rank, it 
suffices to prove surjectivity on fibers at all classical  
points $x$. In other words, if $P_x$ denotes the global sections of 
$(x^\ast\UCF)^\tau$, we have to prove (i) that for any such $x$ the inclusion 
$P\into P_x$ is an equality and (ii) that $P_x\otimes_A\BOne_{x,A(1)}\longto x^\ast\UCF$ 
is an isomorphism. By our assumption on ranks and part~(b) of the 
previous lemma, the ranks of $P$, of $P_x$ and of $x^\ast\UCF$ 
agree. Part~(c) of the lemma then establishes (ii). Since $P$ is 
saturated in $P_x$ by Corollary~\ref{SatUModTau}, having 
the same ranks implies~(i). 
\end{proof} 
 
\medskip 
 
\begin{remark} 
If $\UCF$ arises by coefficient change from a $\tau$-sheaf 
$(\CG,\tau_\CG)$ over $A(\infty)$, we have seen in 
Proposition~\ref{TauInvOnAInfty} that $\UCF^\tau(X) \cong 
\UCG^\tau(X)$. So in particular $\rk_A \UCF^\tau(X) = \rk_A 
\UCG^\tau(X)$. However, it is in general not true that 
\[ 
(\CG,\tau_\CG) \es \cong \es \BOne_{X,A(\infty)} \otimes_A \UCG^\tau(X) 
\] 
even if $\rk_A \UCG^\tau(X) = \rk \CG$. The reason for this is that in 
cases we will be interested in later, the homomorphism $\tau_\CG$ is 
typically not an isomorphism over all of $\FA(\infty)$. Then the support $D$ of the coherent sheaf 
\[ 
\CG / \bigl(\BOne_{X,A(\infty)}\otimes_A \UCG^\tau(X)\bigr) 
\] 
is discrete in $X\times_L \mathfrak{A}(\infty)$ and satisfies $\sigma_{X,\mathfrak{A}(\infty)}^\ast(D) \subset D$. Hence $D$ has a limit point at $\infty$. On the other hand since $\mathfrak{A}(1)$ is compact there may be no limit points of $D$ on $X\times_L \mathfrak{A}(1)$. So $\tau_\CF$ must be an isomorphism over $\FA(1)$. 
\end{remark}


\section{Analytic Triviality of $\tau$-Sheaves} \label{SectAnalytTriv} 
\setcounter{equation}{0}

We now study analytic triviality which is the locally {\'e}tale generalisation of 
triviality of the previous section, and which in fact we define `on {\'e}tale stalks'.  
After describing some basic properties, 
we will show that the locus of analytic triviality in  
$\CM(X)$ is (Berkovich) open and give a first simple example.  
Using the openness property we show that any analytically trivial $\tau$-sheaf  
will be trivial over some temperate {\'e}tale covering (Def.~\ref{DefTemperateFG}).  
 
\medskip 
 
Let $\UCF$ be a locally free rigid analytic $\tau$-sheaf over 
$A(1)$ on $X$. In order to take a more general approach to 
triviality of a $\tau$-sheaf we will consider $\UCF$ over the {\'e}tale 
site. Generalising \cite[(2.3)]{anderson} and motivated by 
Corollary~\ref{CritUnif}, we define: 
\begin{Def} \label{AnTrivDef}
The $\tau$-sheaf $\UCF$ is called {\em analytically trivial} if for every
{\'e}tale analytic point $y\in\CM_\et(X)$ (Sect.~\ref{EtaleSheaves}) the rank of
$(y^*\UCF)^\tau(k(y))$ as an $A$-module is equal to the rank of the
sheaf $y^*\CF$. 
\end{Def} 
 
Clearly Definition~\ref{AnTrivDef} is in accordance with the usual principle that
properties of {\'e}tale sheaves can be defined and established on stalks at
{\'e}tale points. Note however that for general sheaves on rigid analytic spaces
this principle is {\em false\/}. Theorem~\ref{PointwiseTau} below will show
that indeed the above is a sensible definition.


\smallskip 
 
We begin with two simple consequences of the above definition.

\begin{prop}\label{AnalytTrivUnderPull} 
If $\UCF$ is analytically trivial, then so is $\pi^*\UCF$ for every 
general morphism $\pi:X'\to X$. Moreover the $\tau$-sheaves $\UCF^{\otimes n}$, 
$\Sym^k\UCF$ and $\bigwedge^n\UCF$ are analytically trivial. 
\end{prop} 
We omit the obvious proof. 
 
\begin{prop}\label{ExtOfAnTriv} 
Suppose  $0\to\UCF\to\UCG\to \UCH\to 0$ is a short exact sequence of 
locally free rigid analytic $\tau$-sheaves. Then $\UCG$ is 
analytically trivial if and only if the same holds for $\UCF$ 
and~$\UCH$. 
\end{prop} 
\begin{proof} 
We only need to prove the assertion for $X=\Spm F$, where $F$ is a 
complete normed algebraically closed field. This is shown in 
\cite[Lem.~2.7.2]{anderson}. 
\end{proof} 
 
\medskip 
 
Next we  investigate the analytically trivial locus of some $\UCF$ 
on~$X$. 
\begin{lemma}\label{LemOnOpen} 
Suppose $\UCF$ is a rigid analytic locally free $\tau$-sheaf of rank 
$r$ over $A(1)$ on $X$ such that $x^\ast\UCF$ is analytically trivial for 
some $x\in\CM(X)$. Then there exists a wide open neighborhood $U$ of 
$x$ in $X$ such that $\UCF$ is analytically trivial on~$U$. 
\end{lemma} 
Note that if for some $x\in \CM(X)$ there exists a general morphism from a point $y$ to $x$ at which $y^\ast\UCF$ is analytically trivial, then by the 
definition of analytic triviality $x^\ast\UCF$ is analytically trivial. 
\begin{proof} 
As usual, we may assume $A=\BF_q[t]$. 
By our hypothesis there exists a finite separable field extension $L'$ 
of $k(x)$ and elements $\ol{f}_1, \ldots, \ol{f}_r \in (x^\ast\UCF)^\tau(L')$, which form an $\Fz_q[t]$-basis of 
$(x^\ast\UCF)^\tau$. So by 
Lemma~\ref{ExtendOnV}, there exist an 
{\'e}tale morphism $\pi\!:Y\to X$, a point $y\in Y$, a wide open neighborhood 
$V$ of $y$ such that $L'\cong k(y)$, the sections $\ol{f}_j$ extend to 
$V$ and $x=\pi(y)$. By Corollary~\ref{SatUModTau} the rank of 
$(\pi^*\UCF)^\tau(V)$ is again $r$, and so by 
Corollary~\ref{CritUnif}, $\pi^*\UCF$ is analytically trivial on 
$V$. Moreover $U:=\pi(V)$ is a wide open neighborhood of $x$, and it 
follows that $\UCF$ is analytically trivial on~$U$.  
\end{proof} 
 
\begin{thm} \label{UnivOpen} 
Let $X$ be an affinoid rigid analytic space and $\UCF$ a locally free rigid analytic $\tau$-sheaf over $A(1)$ on $X$. Then the subset 
\[ 
M := \{ x \in \CM(X): x^\ast\UCF \enspace \text{is analytically trivial} \} \enspace \subset \enspace \CM(X) 
\] 
is open. 
\end{thm} 
 
\begin{proof} 
This is immediate from the previous corollary. By it we have for 
each $x\in M$ a wide neighborhood $U$ of $x$ on which $\UCF$ 
is analytically trivial. This however means that $y^\ast\UCF$ is 
analytically trivial at all $y\in \CM(U)$, and the latter set is open 
in $\CM(X)$. 
\end{proof} 
 
\begin{remark} 
By the same argument Theorem~\ref{UnivOpen} can be proved more generally for paracompact rigid analytic spaces $X$; cf.~the remark at the end of Section~\ref{AnalyticPoints}. 
\end{remark}

Example~\ref{smallexample} and the example of Section~\ref{SectExample} below raise 
the following question which goes back to a conversation with Pink: 
 
\begin{question} \label{OpenQuestion1} 
Let the notation be as in Theorem~\ref{UnivOpen}. Is the set 
$\CM(X)\smallsetminus M$ a special subset of~$X$? 
\end{question} 
 
If the answer is in the affirmative, this would have strong 
consequences. For instance, since any non-empty special subset contains a 
classical point, a $\tau$-sheaf would be analytically trivial, if it is 
so at all classical points.

\smallskip 
 
With respect to the canonical topology on the affinoid $L$-space $X$ the map $X \to \CM(X)$ is a homeomorphism onto a dense subset. Thus we see: 
 
\begin{cor} \label{UnivOpenCor} 
Let $X$ be an affinoid rigid analytic space and $\UCF$ a locally free rigid analytic $\tau$-sheaf over $A(1)$ on $X$. Then the subset 
\[ 
U := \{ x \in X: x^\ast\UCF \enspace \text{is analytically trivial} \} \enspace \subset \enspace X 
\] 
is open with respect to the canonical topology on $X$. 
\end{cor}

To illustrate this situation we give a simple example. A more involved one will be computed in Section~\ref{SectExample}. 
 
\begin{Bigexample}\label{smallexample} 
Let $A=\Fz_q[t]$, $B=L\langle a,b\rangle$ and $X=\Spm B$. Consider the rigid analytic $\tau$-sheaf $\UCF$ over $A(1)$ on $X$ given by the $B\langle t\rangle$-module $F=B\langle t\rangle$ and $\tau= (a+bt)\cdot\sigma$. Then the subset of $\CM(X)$ over which $\UCF$ is analytically trivial is 
\[ 
\{\,x\in \CM(X): |b(x)|<|a(x)|\,\}\,. 
\] 
 
\begin{proof} 
Fix an analytic point $x\in \CM(X)$. 
We want to find a nonzero $\tau$-invariant $f$ which we write as $f=\sum_{n\geq 0} u_n t^n$, i.e. the coefficients have to satisfy 
\[ 
u_n - a(x)\,u_n^q = b(x)\, u_{n-1}^q \qquad \text{for all}\quad n\geq 0\,. 
\] 
If $a(x) = 0$ then $u_n=0$ for all $n$ and $\UCF$ is not analytically trivial at $x$. 
 
So let $a(x) \neq 0$ and fix a $(q-1)^{\rm st}$ root $\alpha \in k(x)^{alg}$ of $a(x)$. Then setting $v_n=\alpha u_n$ the condition on the coefficients becomes 
\begin{equation} \label{Example1Eq} 
v_n -v_n^q = c \,v_{n-1}^q\,, \qquad c:=\frac{b(x)}{a(x)}\,. 
\end{equation} 
We may assume that $v_0\in \Fz_q\mal$. 
 
If $|b(x)|\geq |a(x)|$ then by induction $|v_n| = |c|^{n/q}\geq 1$ and $\UCF$ is not analytically trivial at $x$. 
 
If $|b(x)|<|a(x)|$ then we may chose $|v_n|=|c\,v_{n-1}^q|$ and by induction $|v_n| = |c|^{1+q+\ldots+q^{n-1}}<1$. Hence $\sum_n v_n t^n \in k(x)\langle t\rangle$ and $\UCF$ is analytically trivial at $x$. This finishes the proof. 
\end{proof} 
 
Let $U$ be the admissible subset of $X$ on which $|b|<|a|$ and let 
$\pi\!:Y\to U$ be the finite {\'e}tale Galois covering given by adjoining a $(q-1)^{\rm st}$ root of $a$. Then $\rk_A \UCF^\tau(Y) = 1$, since we may explicitly solve equation (\ref{Example1Eq}) globally on $Y$. Thus we have 
$\pi^*\UCF \cong \BOne_{Y,A(1)} \otimes_A (\pi^*\UCF)^\tau(Y)$. 
\end{Bigexample}

\bigskip 
 
Let us now give some equivalent conditions for analytic triviality: 
\begin{thm} \label{PointwiseTau} 
Let $X$ be an affinoid rigid analytic space and $\UCF$ a locally free rigid analytic $\tau$-sheaf over $A(1)$ on $X$. Then the following are equivalent: 
\begin{enumerate} 
\item 
$\UCF$ is analytically trivial, 
\item \label{PartBprev} 
for every analytic point $x\in \CM(X)$ the $\tau$-sheaf $x^\ast\UCF$ is analytically trivial, 
\item There exists a finite cover $U_i$ of $X$ by affinoids and 
  {\'e}tale morphisms $\pi_i\!:V_i\to X$ with $\pi_i(V_i)=U_i$, such 
  that the $\UCF_i:=\pi^*_i\UCF$ satisfy \ 
  $(\UCF_i)^\tau(V_i)\otimes_A\BOne_{V_i,A(1)}\cong \UCF_i$. 
\item  
there exists a connected temperate {\'e}tale Galois covering $\pi:Z\to X$ with  
\[ 
\pi^\ast\UCF \enspace \cong \enspace \BOne_{Z,A(1)} \otimes_A \UCF^\tau(Z)\,. 
\] 
\end{enumerate} 
\end{thm} 
 
We do not know whether one can in fact trivialize an analytically 
trivial $\tau$-sheaf over a finite {\'e}tale covering. The following 
example shows that for general (paracompact) rigid analytic spaces, 
this can (obviously) not be expected. 
\begin{example} 
Let $\UCF$ be the $\tau$-sheaf associated to the universal 
Drinfeld-$A$-module of the analytic moduli space $\FM_\Fn^r$ of 
Drinfeld-$A$-modules of some fixed rank $r$ and with some level 
$\Fn$-structure, cf.~\cite{drinfeld}. Let $\Omega^r$ be the Drinfeld 
upper half space as in Example~\ref{FirstExs} and $\pi\!:\Omega^r\to 
\FM_\Fn^r$ be the canonical covering constructed 
in~\cite{drinfeld}. Then $\pi^*\UCF$ is in fact `trivial' in the sense 
of Corollary~\ref{CritUnif}. In particular $\UCF$ is analytically 
trivial. However it may easily be shown that there does not exist a 
finite {\'e}tale cover of $\FM_\Fn^r$ over which $\UCF$ becomes trivial. 
\end{example}

\begin{proof}[Proof of Theorem~\ref{PointwiseTau}] 
The equivalence (a)$\Leftrightarrow$(b) is clear from the definition 
and Corollary~\ref{StalkOfFtauAndOverC}~\ref{StalkOfFtauAndOverC_Parta} and~\ref{StalkOfFtauAndOverC_Partb}. The implication 
(c)$\Rightarrow$(b) is clear from the definitions, and 
(d)$\Rightarrow$(c) follows from the fact that $Z\to X$ is a covering 
for the {\'e}tale topology. We now prove (b)$\Rightarrow$(d). 
 
By Remark~\ref{Taguchi's Trick} we may assume that $A=\BF_q[t]$. 
As in the proof of Lemma~\ref{LemOnOpen}, for $x\in \CM(X)$  there 
exist an {\'e}tale morphism $\pi_x\!:Y_x\to X$, a point $y\in Y_x$ with 
$\pi_x(y)=x$ and  a wide open neighborhood $V_x$ of $y$ in $Y_x$ such that  
$\pi_x^*\UCF|_{V_x}\cong 
(\pi_x^*\UCF)^\tau(V_x)\otimes_A\BOne_{V_x,A(1)}$. Moreover by 
Lemma~\ref{ExtendOnV}~(b), we may assume that $\pi_x$ is the morphism  
$X_{\UCF,t^{N_x}} \to X$ for some $N_x\in\BN$.

By Proposition~\ref{AdmCovering} the $\pi(V_x)$ form an admissible covering of $X$. Since $X$ is quasi-compact, there is a finite set $S$ of analytic points $x$ such that the $\pi(V_x)$ for $x\in S$ form an admissible covering of $X$. Let  
\[ 
N \es = \es \max\{N_x: x\in S\}\,. 
\] 
We set $Y=X_{\UCF,t^N}$. Then $Y\to X$ is a torsor for the group $G=\GL_r(\BF_q[t]/t^{N+1})$. We replace $V_x$ by its preimage under the finite {\'e}tale morphism $Y\to X_{\UCF,t^{N_x}}$. Consider the admissible covering 
\[ 
\CV \es := \es \{gV_x: x\in S, g\in G\} 
\] 
of $Y$. On $gV_x$ we have $\UCF^\tau(gV_x) \cong (g^{-1})^\ast\UCF^\tau(V_x) \cong \UCF^\tau(V_x) \cong \BF_q[t]^r$. These isomorphisms form a \v{C}ech-cocycle 
\[ 
\alpha \es \in \es \CKoh^1\bigl(\CV,\GL_r(\BF_q[t])\bigr)\,. 
\] 
In general one can not expect to trivialize this cocycle by a finite {\'e}tale covering of $Y$. However, it can be trivialized by a topological covering $Z$ of $Y$ which is Galois over $X$. We will construct $Z$. Let $N(\CV)$ be the nerve of the covering $\CV$ of $Y$. Recall that this is a simplicial set whose $n$-simplices correspond to the connected components of the intersections of $n$ distinct elements of $\CV$. There exists an equivalence of categories 
\[ 
\{\text{topological covering spaces of }Y\text{ split by }\CV\} \es\cong\es \{\text{covering spaces of }N(\CV)\}\,. 
\] 
Now let $\wt{N}$ be the universal covering of $N(\CV)$ and let $Z$ be the corresponding topological covering space of $Y$. Clearly $Z$ is Galois over $Y$ with Galois group $\Aut(Z/Y) = \Aut\bigl(\wt{N}/N(\CV)\bigr)$. Furthermore, $Z\to X$ is a temperate {\'e}tale Galois covering. Indeed, as $Z$ is universal for topological covering spaces of $Y$ split over $\CV$, every $g\in G=\Aut(Y/X)$ lifts to an automorphism of $Z$. Therefore we have an exact sequence of discrete groups 
\[ 
1\to \Aut(Z/Y) \to \Aut(Z/X) \to \Aut(Y/X) \to 1\,. 
\] 
From $Z /\Aut(Z/Y) = Y$ and $Y/\Aut(Y/X)=X$ we conclude that $Z/\Aut(Z/X)=X$.  
 
Moreover, the cocycle $\alpha$ is trivialized by $Z$ and hence $\UCF^\tau(Z)\cong \BF_q[t]^r$. (We may assume that $Z$ is connected by replacing $Y$ by one of its connected components.) Now the assertion follows with Prop.\ \ref{CritUnif}. 
\end{proof} 
 
Using for instance the characterization in part~\ref{PartBprev}, the following is immediate, since it holds point-wise:  
\begin{cor} 
If $\UCF$ is analytically trivial, then $\tau_\CF$ is an isomorphism. 
\end{cor} 
 

\section{Exactness of the Fundamental Sequence} \label{GenArtSch} 
\setcounter{equation}{0}

In this section, we give some results on the right exactness of the  
Artin-Schreier type sequence (\ref{ASlikeSeq})  for $A$-coefficients 
from the introduction. Consider Example~\ref{smallexample} and let $x\in \CM(X)$ be an analytic point. The reader may easily verify that the sequence (\ref{ASlikeSeq}) for the fiber $x^\ast\UCF$ over $x$ is right exact either if $x^\ast\UCF$ is analytically trivial or if $\tau=0$. We will shortly see that this is a typical phenomenon. 
 
\smallskip 
 
\begin{prop}\label{AStheoryI} 
Suppose $\UCF$ is analytically trivial over $A(1)$ on $X$. Then the 
sequence $0\to\UCF^\tau\to W(\CF)\stackrel{1-\tau}\to W(\CF)\to0$ of 
{\'e}tale sheaves is exact. 
\end{prop} 
\begin{proof} 
By Corollary~\ref{StalkOfFtauAndOverC}~\ref{StalkOfFtauAndOverC_Partd}, it suffices to prove the 
assertion on fibers at  
{\'e}tale analytic 
points. So let $F$ be a 
complete normed algebraically closed field extension of $L$. As $\UCF^\tau$ is defined to be the kernel of $\id-\tau$, the only part that needs proof is 
the surjectivity of $(1-\tau)$. It clearly suffices to consider the 
case where $\UCF=\BOne_{X,A}$.  
 
Let $f=\sum_{n=0}^\infty b_n t^n\in F\langle t\rangle$ be a section of 
$\CO_{X\times_L\FA(1)}$. The equation $(1-\sigma_{X,\FA(1)})g=f$ leads to the equation \[ 
(\id -\sigma) \bigl(\sum_{n=0}^\infty c_n t^n\bigr) \enspace = \enspace \sum_{n=0}^\infty b_n t^n \,. 
\] 
where $g=\sum_{n=0}^\infty c_n t^n$. Since $F$ is algebraically closed, 
the resulting Artin-Schreier equations $c_n^q-c_n=b_n$ can all be 
solved inside $F$, and the condition $|b_n|\to 0$ for $n\to\infty$, 
implies the same for $|c_n|$. Hence $g$ lies in $F\langle t\rangle$, 
and the surjectivity is established. 
\end{proof} 
 
Our next aim is to generalize the above. We consider $\tau$-sheaves over $A(1)$ on $X$. 
\begin{lemma}\label{AnTrivByNilIsSplit} 
Any extension of an analytically trivial by a locally free nilpotent 
$\tau$-sheaf splits. 
\end{lemma} 
\begin{proof} 
Consider a short exact sequence  
$$0\longto\UCF\longto\UCE\longto\UCG\longto0$$ 
in which $\UCF$ is nilpotent, say $\tau_\CF^n=0$, and locally free, 
and in which $\UCG$ is analytically trivial. Because all the above 
sheaves are therefore locally free, the top sequence in the following 
commutative diagram is exact: 
$$\xymatrix{ 
0\ar[r]&\ar[d]^0(\sigma_{X,\FA(1)}^n)^\ast\CF\ar[r]& 
\ar[d]^{\tau_\CE^n}(\sigma_{X,\FA(1)}^n)^\ast\CE\ar[r]& 
\ar[d]_{\cong}^{\tau_\CG^n}(\sigma_{X,\FA(1)}^n)^\ast\CG\ar@{-->}[dl]\ar[r]&0\\ 
0\ar[r]&\CF\ar[r]&\CE\ar[r]&\CG\ar[r]&0\rlap{.}\\ 
}$$ 
It follows that there is a unique dashed homomorphism. Its composite 
with the inverse of the isomorphism $\tau_\CG^n$ is easily seen to be 
a splitting. 
\end{proof} 
\begin{cor}\label{RepeatedExt} If $\UCF$ is a rigid analytic $\tau$-sheaf over $A(1)$ on $X$ with a 
  filtration such that all subquotients are locally free and either 
  nilpotent or analytically trivial. Then $\UCF$ is an extension of 
  a locally free nilpotent $\tau$-sheaf by an analytically trivial 
  $\tau$-sheaf.  
\end{cor} 
\begin{proof} 
This is immediate from Proposition~\ref{ExtOfAnTriv} and the previous 
lemma. 
\end{proof} 
 
\begin{thm}\label{AStheoryII} 
Let $X$ be a rigid analytic space and 
$\UCF$ a locally 
free rigid analytic $\tau$-sheaf over $A(1)$ on $X$. Assume that for 
every analytic point $x$ of $X$ the $\tau$-sheaf $x^\ast\UCF$ is 
an extension of a locally free nilpotent by an analytically trivial 
$\tau$-sheaf. Then 
\[ 
0 \longrightarrow \UCF^\tau \longrightarrow W(\CF) 
\xrightarrow{\id -\tau_\CF} W(\CF) \longrightarrow 0 
\] 
is an exact sequence of {\'e}tale sheaves on $X$. 
\end{thm} 
 
\begin{proof} 
As $\UCF^\tau$ is defined to be the kernel of $\id -\tau_\CF$, we only 
have to show that $\id -\tau_\CF$ is surjective. By 
Corollary~\ref{StalkOfFtauAndOverC}~\ref{StalkOfFtauAndOverC_Partd}, it suffices to prove 
surjectivity on fibers at 
{\'e}tale analytic 
points. Suppose therefore that $X=\Spm F$ for some complete normed 
algebraically closed field. 
 
By assumption $\UCF$ is an extension 
of a nilpotent locally free by an analytically trivial sheaf. By an 
obvious $9$-term square, it suffices to prove surjectivity separately 
for the latter two cases. If $\UCF$ is nilpotent, surjectivity is 
clear. In the other case, it was shown in Proposition~\ref{AStheoryI}.  
\end{proof} 
 
\begin{question} 
We pose the following problem which in spirit is similar to  
Question~\ref{OpenQuestion1} and the remark following it: 
Suppose one assumes 
in the above theorem that $x^\ast\UCF$ is an extension of a 
locally free nilpotent by an analytically trivial $\tau$-sheaf 
only at every classical point $x$ of $X$. Does the conclusion of the 
theorem still hold? 
 
Another interesting question is to give conditions under which one  
can recover $\UCF$ from the locally {\'e}tale sections of $\UCF^\tau$  
-- at least up to nilpotent parts of~$\UCF$. 
\end{question}


\section{An Example} \label{SectExample} 
\setcounter{equation}{0}

The following example is due to R. Pink. 
 
Let $\zeta \in K_\infty$ satisfy $|\zeta|<1$ and let  
\[ 
X = \Spec K_\infty[a,b,c,d]/(a+d +2\zeta\,,\, ad-bc-\zeta^2)\,. 
\] 
We set  
\[ 
\Delta_1 =\left(\begin{array}{cc} a& b\\ c&d \end{array}\right) 
\] 
 and $\Delta = \Id + t\Delta_1$. Then $\det \Delta = (1-\zeta t)^2$. Let $X^{\rig}$ be the analytification of $X$. As in the remark at the end of Section~\ref{AnalyticPoints} we consider analytic points $x$ of $X^{\rig}$ and we represent them by their corresponding matrices $\Delta(x)$ and $\Delta_1(x)$. 
 
On $X$ we define a (locally) free algebraic $\tau$-sheaf $\UCF$ over $A=\Fz_q[t]$  by  
\[ 
\CF = \CO_{X\times_{\Fz_q} \Spec A}^{\oplus 2}\,,\quad \tau = \Delta \cdot \sigma\,. 
\] 
 
There is an action of $\GL_2(\Fz_q)$ on $X$ and $X^{\rig}$ given by  
\[ 
\GL_2(\Fz_q)\,\times\,X \longrightarrow X\,,\quad (g,\Delta) \mapsto g\Delta g^{-1}\,. 
\] 
Under this action we have isomorphisms $g^\ast\UCF\cong\UCF$ for all $g\in\GL_2(\BF_q)$. We denote the associated rigid analytic $\tau$-sheaf over $A(1)$ by $\wt\UCF$. Note that $\UCF$ is the $\tau$-sheaf associated to the family of $\BF_q[t]$-modules $(E,\phi)$ of dimension $2$ and rank $2$ with $E=\CO_X^2$ and 
\[ 
\phi_t = -\Delta_1^{-1} + \Delta_1^{-1}\sigma\,, 
\] 
when we take $\zeta^{-1}$ as the image of $t$ in $K_\infty$; cf.~Section~\ref{SectFamAMod}. Moreover $X$ together with $\UCF$ is a fine moduli space for ``polarized'' $t$-motives of rank $2$ and dimension $2$ with appropriate level structure; cf.~\cite[{\S} 4]{hartl1}.  
 
The following result illustrates Theorem~\ref{UnivOpen}. 
 
\begin{prop} \label{ExampleProp} 
\begin{eqnarray*} 
&&\{\,y\in \CM(X^{\rig}): y^\ast\wt\UCF \enspace\text{\rm is analytically trivial} \,\} \\[0.2cm] 
&& \qquad\qquad\qquad\enspace = \enspace \bigcup_{g\in\GL_2(\Fz_q)} g\,\{\,y \in \CM(X^{\rig}) : \;|a(y)|,|c(y)|,|d(y)| <1\,\}\,. 
\end{eqnarray*} 
\end{prop} 

\medskip

Another explicit family of $t$-motives for Theorem~\ref{UnivOpen} had been
given by F.\ Gardeyn \cite[II.2]{GardeynDiss} in his thesis. His example can
be considered as a special case of the above.

\begin{proof} 
We investigate the $\tau$-invariants of $\wt\UCF$ which we write as $2\times 2$-matrices $\Phi = \Delta \,{}^{\sigma\!} \Phi$ with $\Phi = \sum_{n\geq 0} \Phi_n t^n$. We set $\Phi_n = \left(\begin{array}{cc} u_n & v_n \\ w_n & x_n\end{array}\right)$ and $\Phi_0=\Id$. This implies for $n\geq 1$ 
\begin{eqnarray} \label{ExampleEq1} 
\nonumber \Phi_n -{}^{\sigma\!} \Phi_n & = & \Delta_1 \,{}^{\sigma\!} \Phi_{n-1} \qquad \text{or equivalently}\\[0.2cm] 
u_n - u_n^q & = & a\, u_{n-1}^q + b\, w_{n-1}^q \,,\\[0.2cm] 
\nonumber w_n - w_n^q & = & c\, u_{n-1}^q + d\, w_{n-1}^q \,,\\[0.2cm] 
\nonumber v_n - v_n^q & = & a\, v_{n-1}^q + b\, x_{n-1}^q \,,\\[0.2cm] 
\nonumber x_n - x_n^q & = & c\, v_{n-1}^q + d\, x_{n-1}^q \,. 
\end{eqnarray} 
The $\tau$-sheaf $y^\ast\wt\UCF$ is analytically trivial if and only if we can find solutions $\Phi_n(y)$ in $M_2\bigl(k(y)^\alg\bigr)$ to these equations with $|\Phi_n(y)| \to 0$ for $n\to \infty$. 
 Fix a point $y$ of $\CM(X^\rig)$. To ease notation we write $a$ instead of $a(y)$, etc. 

We first prove the inclusion ``$\supset$''.   
If $|\Delta_1|<1$ then the product $\Phi := \Delta \cdot \,{}^{\sigma\!}\Delta \cdot \,{}^{\sigma^2\!}\Delta \cdot \ldots\;$ converges in $\GL_2\bigl(k(y)\langle t\rangle\bigr)$ and hence $\wt\UCF$ is analytically trivial at $y$.  
 
Next observe that $\wt\UCF$ is analytically trivial at $y$ if and only if it is at $g(y)$ for every $g \in \GL_2(\Fz_q)$. So it suffices to treat the following case.
 
Let $r\in K_\infty^{alg}$ be a constant with $|\zeta|\leq |r|<1$ such that $|a|, |d| \leq |r|$ and $|r|\le|b|$. This implies in particular  
\[ 
|c| \es= \es\frac{|bc|}{|b|} \es=\es \frac{|ad -\zeta^2|}{|b|} \es\leq\es
 \frac{|r|^2}{|b|\;\,}\es\le\es |r|\es <\es 1\,. 
\] 
Suppose we choose the solutions of (\ref{ExampleEq1}) in such a way that
\[
\begin{array}{c@{\;<1\quad\Rightarrow\quad}c @{\;=\;}c}
|u_n-u_n^q| & |u_n| & |u_n-u_n^q|\,,\\[1mm]
|v_n-v_n^q| & |v_n| & |v_n-v_n^q|\,,\\[1mm]
|w_n-w_n^q| & |w_n| & |w_n-w_n^q|\,,\\[1mm]
|x_n-x_n^q| & |x_n| & |x_n-x_n^q|\,,
\end{array}
\]
for $n\ge 1$. Note that this is possible, because for $\alpha\in k(y)^\alg$
with $|\alpha|<1$, the equation $x-x^q=\alpha$ has the solution
$x=\sum_{\nu=0}^\infty \alpha^{q^\nu}$. In this situation we have the estimates
\[
\begin{array}{c@{\es\le\es\max\{\;}c@{\;,\;}c@{\;\}\,}l}
|u_n|^q & |a|\,|u_{n-1}|^q & |b|\,|w_{n-1}|^q & , \\[1mm]
|w_n|^q & |c|\,|u_{n-1}|^q & |d|\,|w_{n-1}|^q & , \\[1mm]
|v_n|^q & |a|\,|v_{n-1}|^q & |b|\,|x_{n-1}|^q & , \\[1mm]
|x_n|^q & |c|\,|v_{n-1}|^q & |d|\,|x_{n-1}|^q & . 
\end{array}
\]
These estimates are crude because we are not taking advantage of the gain we
get when the right side is of absolute value less than one. From this we
obtain the even cruder estimates 
\begin{eqnarray*}
|u_n|^q & \le & \max\{\,|r|\,|u_{n-1}|^q\;,\;|b|\,|w_{n-1}|^q\,\}\\[2mm]
& \le & |r|\,\max\{\,|u_{n-1}|^q\;,\;|b/r|\,|w_{n-1}|^q\,\}\,,\\[5mm]
|b/r|\,|w_n|^q & \le & |b/r|\,\max\{\,|r^2/b|\,|u_{n-1}|^q\;,\;|r|\,|w_{n-1}|^q\,\}\\[2mm]
& \le & |r|\,\max\{\,|u_{n-1}|^q\;,\;|b/r|\,|w_{n-1}|^q\,\}\,.
\end{eqnarray*}
These yield \begin{equation}\label{RecursEstimate}
\max\{\,|u_{n}|^q\;,\;|b/r|\,|w_{n}|^q\,\}\le
|r|^n\max\{\,|u_{0}|^q\;,\;|b/r|\,|w_{0}|^q\,\},\end{equation}
 so that $|u_n|\to0$ and $|w_n|\to 0$ as $n\to\infty$. A similar estimate shows that $|v_n|\to0$ and $|x_n|\to0$. We conclude that $y^\ast\wt\UCF$ is analytically trivial.
This finishes the proof of the inclusion ``$\supset$''. (In Remark~\ref{RemRichardsExample} we explain how this example also illustrates Theorem~\ref{PointwiseTau}.)
 
\medskip 
 
For the opposite inclusion fix an analytic point $y$ of $X^{\rig}$ which does not belong to the right hand side of the claimed equality in Proposition~\ref{ExampleProp}. Thus $|a|=|d|\ge1$. We shall show that $y^\ast\wt\UCF$ is not analytically trivial. 
Let us introduce the following new variables
\[
\wt\Delta_1\es=\es
\left(
\begin{array}{cc}
\tilde a & \tilde b \\ \tilde c & \tilde d
\end{array}\right)\es:=\es\left(
\begin{array}{cc}
-\zeta-a & -b \\ -c & -\zeta-d
\end{array}\right)\es=\es-\zeta\Id-\Delta_1 \,.
\]
We have the relations
\[
\wt\Delta_1^2\;=\;0\,,\qquad \tilde b\wt\Delta_1\;=\; \left(
\begin{array}{cc}
-\tilde b\tilde d & \tilde b^2 \\ -\tilde d^2 & \tilde b \tilde d
\end{array}\right)\;=\;\left(
\begin{array}{c}
\tilde b \\ \tilde d
\end{array}\right) \bigl(-\tilde d\quad\tilde b\bigr)\,.
\]
Observe furthermore that $|a|=|d|\ge1$ implies $|\tilde a|=|a|,|\tilde b|=|b|, |\tilde c|=|c|$ and $|\tilde d|=|d|$.

Let us now assume that $y$ is chosen such that the value $|\tilde d(y)|$ is minimal among all $\bigl|\tilde d\bigl(g(y)\bigr)\bigr|$ for $g \in \GL_2(\Fz_q)$ and such that $|\tilde b(y)|\ge|\tilde c(y)|$.  In particular we have $|\tilde b|\ge|\tilde d|=|\tilde a|\ge1$.
Note that the minimality assumption for $|\tilde d|$  implies
\begin{equation} \label{EqBAndD}
|\tilde b^{q-1} - \tilde d^{q-1}| = |\tilde b^{q-1}|\,. 
\end{equation}
Indeed if not then writing $\tilde b^{q-1} - \tilde d^{q-1}\;=\;\prod_{\alpha\in\BF_q\mal}(\tilde b-\alpha\,\tilde d)$ yields an $\alpha \in \Fz_q\mal$ with $|\tilde b- \alpha\,\tilde d|<|\tilde b|$. In particular $|\tilde b|=|\tilde d|$. But then
\[
\left(\begin{array}{cc} 1 & 0 \\ -\alpha^{-1}  & 1\end{array}\right)\,\wt\Delta_1\,
\left(\begin{array}{cc} 1 & 0 \\ \alpha^{-1}  & 1\end{array}\right)\es=\es
\left(\begin{array}{cc} \ast & \ast \\ \ast  & \tilde d-\alpha^{-1}\tilde b\end{array}\right)\,.
\]
This contradicts the minimality of $|\tilde d|$.

\medskip

Next one shows by induction that
\begin{eqnarray*}
\left(\begin{array}{c} v_n^q-v_n \\ x_n^q-x_n \end{array}\right) & = & -\Delta_1\left( \begin{array}{c} v_{n-1}^q \\ x_{n-1}^q \end{array}\right) \\[1mm]
& = & -\Delta_1\left( \begin{array}{c} v_{n-1}^q -v_{n-1}\\ x_{n-1}^q-x_{n-1} \end{array}\right)\;-\;\Delta_1\left( \begin{array}{c} v_{n-1} \\ x_{n-1} \end{array}\right)\\[1mm]
& = & \sum_{j=1}^n(-\Delta_1)^j\left( \begin{array}{c} v_{n-1} \\ x_{n-1} \end{array}\right)\,.
\end{eqnarray*}
We have $(-\Delta_1)^j\,=\,(\zeta+\wt\Delta_1)^j\,=\,\zeta^j+j\zeta^{j-1}\wt\Delta_1$. Setting $y_n := (-\tilde d\quad\tilde b)\left(\begin{array}{c} v_n \\x_n \end{array}\right)$ we can write
\begin{eqnarray} \label{EqVAndX}
\left(\begin{array}{c}
v_n^q-v_n \\ x_n^q-x_n \end{array}\right) & = &
\left(\begin{array}{c} 1 \\ \tilde d/\tilde b\end{array}\right) y_{n-1} +\\
& & 
+\zeta \left(\begin{array}{c} v_{n-1} \\ x_{n-1}\end{array}\right) +\sum_{j=2}^n\Bigl(j\zeta^{j-1}\left(\begin{array}{c} 1 \\ \tilde d/\tilde b\end{array}\right) y_{n-j}
+\zeta^j\left(\begin{array}{c} v_{n-j} \\ x_{n-j}\end{array}\right)\Bigr)\nonumber
\end{eqnarray}
We will see that the terms in the second line have smaller absolute values
than the corresponding term in the first line. More precisely we make the following
\begin{tabbing}
Claim:\quad for all $n\ge 0$ we have\qquad \= $|x_n| \es=\es |\tilde d|^{(1-q^{-n})/(q-1)} \es\ge\es |x_{n-1}|\es\ge\es1$\\[1mm]
\>$|v_n| \es\le \es|\tilde b|\,|x_n|$\\[1mm]
\>$|y_n| \es= \es|\tilde b|\,|x_n|$\,.
\end{tabbing}
From the claim it follows that $\sum_n \Phi_n t^n$ does not converge for $|t|=1$ and that therefore $x^\ast\wt\UCF$ is not analytically trivial. 

We prove the claim by induction on $n$.
For $n=0$ we have $x_0=1,v_0=0,y_0=\tilde b$.
So now assume the claim for $n-1$. From (\ref{EqVAndX}) we obtain
\[
|v_n^q-v_n|\es = \es |\tilde b|\,|x_{n-1}|\es\ge\es1\qquad\text{and}\qquad |x_n^q-x_n|\es =\es|\tilde d|\,|x_{n-1}| \es\ge\es1\,.
\]
This implies the claimed estimate for $|x_n|$ and the estimate
\[
|v_n|\es = \es \bigl|\tilde b/\tilde d\bigr|^{1/q}\,|x_n|\es\le\es|\tilde b|\,|x_n|\,.
\]
It remains to determine the absolute value of $y_n=\tilde
b\,x_n-\tilde d\,v_n$. Multiplying (\ref{EqVAndX}) with $(-\tilde d\quad\tilde b)$ from the left yields 
$\tilde b \,x_n^q-\tilde d \,v_n^q=\sum_{j=0}^n \zeta^jy_{n-j}$. Together with $y_n^q=\tilde
b^q\,x_n^q-\tilde d^q\,v_n^q$ we obtain:
\[
\left( \begin{array}{cc} -\tilde d^q & \tilde b^q \\ -\tilde d & \tilde b
\end{array}\right)\cdot
\left(  \begin{array}{c} v_n^q \\ x_n^q \end{array}\right)=
\left( \begin{array}{c} y_n^q \\ \sum_{j=0}^n \zeta^jy_{n-j}
\end{array}\right)\,.
\]
By (\ref{EqBAndD}) the determinant of the $2\times2$-Matrix is of absolute value $|\tilde b|^q\,|\tilde d|$. Solving for
$v_n^q$, we find:
$$|y_{n-1}|=|\tilde b|\,|x_{n-1}|=|v_n^q|\le\max\Big\{
\frac{|\tilde b|}{|\tilde d|}\,\Big(\frac{|y_n|}{|\tilde
  b|}\Big)^q\,,\frac{|\,y_{n}|}{|\tilde
  d|}\,,\frac{|\zeta|^1|\,y_{n-1}|}{|\tilde d|}\,,\ldots,
\frac{|\zeta|^n|\,y_{0}|}{|\tilde d|}\Big\}\,.$$
As $|\tilde d|\ge1$, $|\zeta|<1$, and $|y_j|$ is increasing for
$j=0,\ldots,n-1$, we either have
$$ |y_n|\ge |x_{n-1}|\,|\tilde b|\,|\tilde d| \quad\hbox{or}\quad
\frac{|y_n|}{|\tilde b|} \ge (|\tilde d|\,|x_{n-1}|)^{1/q}=|x_n|.$$
Since $|\tilde d|\ge1$, the formula $|x_n|=|\tilde d|^{(1-q^{-n})/(q-1)}$
gives $|x_{n-1}|\,|\tilde d|\ge |x_n|$, so that we have shown $|y_n|\ge
|x_n|\,|\tilde b|$. The converse inequality is clear from 
$$|y_n|=|\tilde b\,x_n-\tilde d\,v_n|\le \max\{|\tilde b|\,|x_n|,
|\tilde d|\,|v_n|\}\le \max\{|\tilde b|\,|x_n|,\bigl|\tilde b/\tilde d\bigr|^{1/q}\,
|\tilde d|\,|x_n|\}= |\tilde b|\,|x_n|.$$
This concludes the proof of Proposition~\ref{ExampleProp}. 
\end{proof} 
 
\medskip

\begin{remark}\label{RemRichardsExample}
We want to explain the relation to Theorem~\ref{PointwiseTau}. Consider constants $r,s \in K_\infty^{alg}$ with $|\zeta|\leq |r|<1<|s|$ and the affinoid subdomain of $X^\rig$
\[
\Spm B\es=\es\{\,y\in X^{\rig}:\;|a(y)|, |c(y)|, |d(y)| \leq |r|\,,\;|r|\leq |b(y)| \leq |s|\,\}\,.
\]
We have just shown that for every analytic point $y\in \CM(B)$ the $\tau$-sheaf $y^\ast\wt\UCF$ is analytically trivial. 

Let $N\in \BZ$ be minimal with $|r^N s|<1$ and define the finite {\'e}tale $B$-algebra $B'$ as 
\[ 
B \,\Big\langle \,\frac{v_1}{s^{1/q}}\,,\, \frac{v_2}{(rs)^{1/q}}\,, \ldots,\, \frac{v_N}{(r^{N-1}s)^{1/q}} \,\Big\rangle 
\] 
modulo the relations 
\[ 
v_1-v_1^q=b\,,\qquad v_n-v_n^q= a\,v_{n-1}^q + b\,x_{n-1}^q\,,\qquad x_n = \sum_{\nu=0}^\infty \bigl(c\, v_{n-1}^q+d\,x_{n-1}^q\bigr)^{q^\nu}\,. 
\] 
Note that the series defining $x_n$ converges in $B'$ by the estimates below. 
Namely, the estimates from (\ref{RecursEstimate}) when applying the initial
values for $\Phi_0$ (and analogous estimates for the $|x_n|$ and $|v_n|$)
imply:
$$\max\{\,|u_{n}|^q\;,\;|b/r|\,|w_{n}|^q\,\}\le
|r|^n\quad \hbox{and}\quad
\max\{\,|v_{n}|^q\;,\;|b/r|\,|x_{n}|^q\,\}\le
|r|^{n-1} \,|b|\,.$$
Thus $|u_{n}|\,,\,|w_{n}|\,,\,|x_n|\,\le\, |r|^{n/q}$ and $|v_{n}|\,\le
(|s|\,|r|^{(n-1)})^{1/q}$. Because $|d|,\,|c|\,<\,1$, a simple inductive
argument proves the convergence of the series for $x_n$. 
\forget{
working the crude estimates we used in the proof of Proposition~\ref{ExampleProp} in their optimal form one can show that there are solutions in $B'$ to the equations (\ref{ExampleEq1}) which for all $z\in \Spm B'$ and all $n\geq 1$ satisfy 
\begin{eqnarray*} 
|u_n(z)| &\leq& |r|^{1+q+ \ldots +q^{n-1}}\,, \\[0.1cm] 
|w_n(z)| &\leq& |r|^{1+q+ \ldots +q^{n-1}}\,\frac{|r|}{|b(z)|}\,, \\[0.1cm] 
|v_n(z)| &\leq& |r^{n-1}b(z)|^{1/q}\enspace \leq \enspace |r^{n-1}s|^{1/q}\,, \\[0.1cm] 
|x_n(z)| &\leq& |r|^n\,. 
\end{eqnarray*} }
From this we also see
that $B'$ is indeed finite over $B$.  Moreover these estimates imply that the series $\Phi = \sum_{n\geq 0} \Phi_n t^n$ converges in $M_2\bigl(B'\langle t\rangle\bigr)$. So the columns of $\Phi$ generate an $A$-submodule of $\wt\UCF^\tau(\Spm B')$ of rank $2$. 
The morphism $\Spm B'\to \Spm B$ is a Galois covering with Galois group 
\[ 
G\es=\es \Bigl\{ \,g = \left(\begin{array}{cc} 1 & \DS\sum_{n=1}^N g_n t^n \\[0.2cm] 0 & 1 \end{array}\right)\Bigr\} \subset \GL_2(\Fz_q[t]) 
\] 
 operating on $B'$ via 
\[ 
g {v_n \choose x_n} \es=\es {v_n\choose x_n} + \sum_{i=1}^n g_i {u_{n-i}\choose w_{n-i}}\,. 
\] 
So in this example the temperate \'etale Galois covering of $\Spm B$ from Theorem~\ref{PointwiseTau} can be chosen finite. However, we do not know whether this is true in general.
\end{remark}

\section{Families of Abelian $A$-Modules} \label{SectFamAMod} 
\setcounter{equation}{0}


\subsection{Definitions} \label{SectFamAMod-Def}

Let $X$ be a rigid analytic space over $L$. We introduce a category of rigid analytic vector group spaces over $X$. Namely let $\bVec(X)$ be the category of locally free coherent $\CO_X$-modules (by the usual abuse of language called vector bundles) with $\Fz_q$-linear {\em algebraic\/} group homomorphisms $\alpha: E \to E'$, i.e., locally with respect to an admissible affinoid covering $\Spm B_i, i\in I$ of $X$ the morphism $\alpha$ can be given by a polynomial 
\[ 
\sum_{\nu=0}^{m_i}A_{i,\nu} \sigma^\nu: B_i^n \longrightarrow B_i^{n'} 
\] 
where $A_{i,\nu} \in M_{n,n'}(B_i)$, the restrictions of $E$ and $E'$ to $\Spm B_i$ are free of rank $n$ and $n'$, respectively, and $\sigma$ denotes the $q$-power map on the components of $B_i^n$. The objects of $\bVec(X)$ are viewed as rigid analytic group spaces over $S$, and therefore the object $\CO_X$ is also called $\Gr_a$. 
 
For an object $E$ of $\bVec(X)$ we let $\Lie E$ be the tangent space of the group $E$ along the zero section. Clearly $\End(\Lie E) \cong \End_{\CO_X}(E)$ and hence there is a diagonal action of $a \in A$ on $\Lie E$ by the inclusion $\iota\!:A \subset K_\infty \subset L \subset \Gamma(X,\CO_X)$. Furthermore, for every morphism $\alpha:E\to E'$ in $\bVec(X)$ we denote its differential by $\alpha':\Lie E \to \Lie E'$. 
 
 
\begin{Def}\label{DefAMod} 
Let $X$ be a rigid analytic space over $L$. An {\em abelian $A$-module $(E,\phi)$ of dimension $d$ and rank $r$ on $X$\/} consists of a vector bundle $E$ of rank $d$ on $X$ and a ring homomorphism 
\[ 
\phi: A \longrightarrow \End_{\bVec(X)}(E): a \mapsto \phi_a\,, 
\] 
such that the following hold: 
\begin{enumerate} 
\item\label{CondAonEMod} 
There exists an admissible covering $\{X_i: i\in I\}$ of $X$ and natural numbers $n_{i,a}$ for each $a\in A$ satisfying 
\[ 
(\iota(a) - \phi_a')^{n_{i,a}} \Lie E|_{X_i} = 0\,. 
\] 
\item 
The sheaf $M(E) := \Hom_{\bVec(X)}(E, \Gr_a)$, equipped with the action of $a \in A$ by composition on the right with $\phi_a$, is a finitely generated $\CO_X \otimes_{\Fz_q} A$-module. (see also Proposition~\ref{LocFreeAmod}.) 
\item\label{CondConEMod} 
$M(E)$ is a flat $\CO_X$-module and the fiber of $M(E)$ at every point $x \in X$ is locally free over $k(x) \otimes_{\Fz_q} A$ of rank $r$. (see also Proposition~\ref{LocFreeAmod}.) 
\end{enumerate} 
\end{Def} 
 
An abelian $A$-module of dimension $1$ and rank $r$ is simply a 
Drinfeld-$A$-module of rank~$r$. 
 
\smallskip 
 
We define an operation $\tau$ on $M(E)$ as composition on the left 
with the Frobenius $\sigma$ on $\Gr_{a,X}$. This makes $M(E)$ into an 
algebraic $\tau$-sheaf over $A$ on $X$. Due to the following proposition the rigid analytic $\tau$-sheaves over $A(1)$ and over $A(\infty)$ associated to $M(E)$ are locally free. 
 
As in Anderson \cite{anderson} one proves: 
\begin{prop}\label{TModActsOnLie} As $\CO_X\otimes A$-modules, the dual of 
  $\Lie E$, and $M(E)/\tau M(E)$ are isomorphic with $A$ acting on 
  $M(E)$ as defined above and on $\Lie E$ via~$\phi'$. 
\end{prop} 
 
\smallskip 
 
The following result is clear from the definition: 
\begin{prop} 
If $\pi:Y\to X$ is a general morphism of rigid analytic spaces, then the pullback of an abelian $A$-module $(E,\phi)$ on $X$ along $\pi$ is an abelian $A$-module on~$Y$. 
\end{prop} 
 
\smallskip 
 
We end this section by giving an equivalent condition to~\ref{CondConEMod}. 
\begin{prop}\label{LocFreeAmod} 
Suppose $\CM$ is a finitely generated $\CO_X \otimes_{\Fz_q}A$-module, where $X=\Spm B$ is an affinoid and $M$ is the $B\otimes_{\BF_q} A$-module corresponding to $\CM$. Then $M$ is locally free over $B\otimes_{\BF_q} A$ of rank $r$ if and only if $M$ is flat over $B$ and for each $x$ in the set of maximal ideals $\Max(B)$ the module $M\otimes_B k(x)$ is locally free over $k(x)\otimes_{\Fz_q}A$ of rank~$r$. 
\end{prop} 
 
\begin{proof} 
Since one direction is obvious, we now assume that $M$ is flat over $B$ and that $M\otimes_B k(x)$ is locally free over $k(x)\otimes A$ for any $x\in\Max(B)$. It suffices to show that for any maximal ideal $\Fn$ of $B\otimes A$ the module $M_\Fn$ is free over $(B\otimes A)_\Fn$ of rank $r$.  
 
Since $B$ is Jacobson and $B\otimes A$ is a finitely generated $B$-algebra $\Fm:=\Fn\cap B$ is a maximal ideal in $B$; cf.~
\cite[Thm.\ 4.19]{Eisenbud}. It therefore suffices to show that for each $\Fm\in\Max(B)$ the module $M_\Fm$ is locally free over $B_\Fm\otimes A$ of rank $r$. Finally because $M$ is finitely generated over $B\otimes A$ and $M\otimes_B k(x)$ is locally free over $k(x)\otimes A$ of rank $r$, the proposition is proved once we show that $\Tor_i^{B_\Fm\otimes A}(M_\Fm,N)=0$ for all finitely generated $B_\Fm\otimes A$-modules $N$ and all~$i>0$.  
 
By the local criterion of flatness, cf.~
\cite{matsu}, it suffices to prove the latter assertion for $k(x)\otimes A$-modules $N$ only. For such there is the following spectral sequence: 
$$\Tor_i^{k(x)\otimes A}(\Tor_j^{B_\Fm}(M_\Fm,k(x)), N)\Longrightarrow \Tor_{i+j}^{B_\Fm\otimes A}(M_\Fm,N).$$ 
By flatness of $M$ over $B$, the spectral sequence degenerates, and because 
$M\otimes_B k(x)$ is flat over $k(x)\otimes A$, except for $i=j=0$ all terms vanish. 
\end{proof}

\subsection{Torsion Points of Abelian $A$-Modules}

\begin{Def} 
The {\em {\'e}tale sheaf $E[\Fa]$ on $X$ of $\Fa$-torsion points of $(E,\phi)$\/} 
is the {\'e}tale sheaf of $A/\Fa$-modules defined for any {\'e}tale 
morphism $Y\to X$ by 
\[ 
E[\Fa] \,(Y)\enspace :=\enspace \bigl\{s\in E(Y):\enspace \phi_a(s)=0 \enspace \text{for every }a\in \Fa\bigr\}\,. 
\] 
\end{Def}

\begin{thm} \label{ComparisonF-E} 
Let $(E,\phi)$ be an abelian $A$-module on $X$, and suppose $\UCF$ is 
the rigid analytic $\tau$-sheaf associated to $M(E)$ as in the first 
construction of Subsection~\ref{Functors}. Then there is a natural 
isomorphism of sheaves of $\nul{A/\Fa}$-modules on~$X$ 
\[ 
E[\Fa] \rbij \UCF[\Fa]\,,\quad s\longmapsto \bigl(h_s:f \mapsto f(s) \bigr)\,. 
\] 
 
Moreover for any ideal $\Fa$, the homomorphism $\tau$ on 
$\UCF/\Fa\UCF$ is an isomorphism, and hence there exists $\GL_r(A/\Fa)$-torsor $\pi\!:Y\to X$ so that 
$\pi^*E[\Fa]$ is isomorphic to $\nul{A/\Fa}^r$. 
\end{thm} 
 
\begin{proof} 
For every {\'e}tale $\pi:Y\to X$ the morphism 
\[ 
E[\Fa](Y) \longto \UCF[\Fa](Y), \quad s\longmapsto \bigl(h_s:f\mapsto f(s)\bigr) \es\in \Hom_{\CO_Y}\bigl(\pi^\ast(\CF/\Fa\CF), \CO_Y\bigr) 
\] 
is an isomorphism. This can be seen locally on $X$ where $E$ is free by choosing an isomorphism $\kappa:E\rbij \CO_X^d$ and applying $h_s$ to the components $\kappa_i$ of $\kappa$ which are elements of $\CF$. 
 
In view of Theorem~\ref{ThmaTorsion}, it remains to prove that $\tau$ is 
bijective on $\UCF/\Fa\UCF$, and it suffices to verify this for principal $\Fa=(a)$.  
For this observe that by condition (a) on an abelian $A$-module the homomorphism  
$\phi'_a$ is bijective on $\Lie E$. By Proposition~\ref{TModActsOnLie}, multiplication by $a$ is  
bijective on $\CF/\tau\CF$, i.e., $\CF=a\CF+\tau\CF$. But then $\tau$ is surjective on $\CF/a\CF$. 
Since $\CF/a\CF$ is locally free the asserted bijectivity follows. 
\end{proof} 
 
\forget{ 
\bigskip 
 
\begin{remark} 
Using that for every $a\in A$ the endomorphism $\phi_a'$ of $\Lie E$ is invertible one can prove that $E[\Fa]$ is representable by a rigid analytic space, which is {\'e}tale over $X$. On the other hand, by comparing the defining equations of $E[\Fa]$ and $\UCF[\Fa]$ one can can proof the stronger result $E[\Fa]\cong \UCF[\Fa]$. For the later we will now give a slightly different proof. 
\end{remark} 
 
}


\subsection{The Exponential Map} \label{SectFamAMod-Exp} 
 
In this section we construct an exponential map for $E$; cf.\ Anderson \cite[{\S} 2.1]{anderson}. As opposed to working over a field, a difficulty in the relative situation lies in the fact that we may not assume that the $\phi_a$ are in Jordan normal form. Therefore we need slightly stronger estimates than Anderson. 
 
\medskip 
 
Let $E$ be an abelian $A$-module of dimension $d$ on $X$. We choose an affinoid refinement of the covering introduced in \ref{CondAonEMod} of Definition~\ref{DefAMod}, such that on each of its open sets $\Spm B$ the sheaf $E|_{\Spm B}$ is free. Let $\kappa: \Gamma(\Spm B, E) \to B^d$ be a fixed isomorphism and $\kappa'$ the induced isomorphism $\Gamma(\Spm B, \Lie E) \to B^d$. We call $\kappa$ and $\kappa'$ {\em coordinates\/} of $E|_{\Spm B}$ and $\Lie E|_{\Spm B}$, respectively. With respect to these coordinates we can represent $\phi_a$ and $\phi_a'$ by matrices 
\begin{equation} 
\phi_a = \sum_{k=0}^{s_a} G_k(a) \sigma^k\,, \qquad \phi_a' = G_0(a)\,,\qquad G_k(a) \in M_d(B)\,, 
\end{equation} 
where $G_0(a) - \iota(a)\Id$ is nilpotent and $\sigma$ acts component-wise as the $q$-power map on $B^d$. We use the abbreviation $z^{(k)} := {}^{\sigma^k}\!z$ for vectors, matrices $z$, etc. 
 
\begin{prop}\label{DefExpMap} 
There exists a unique sequence of matrices $(e_j)$ in $M_d(B)$ with $e_0=\Id$ and $\lim_{j\to \infty} q^{-j} \log |e_j|=-\infty$ defining a rigid analytic function 
\[ 
\exp_{E,\Spm B}: B^d \longrightarrow B^d\,,\quad z \mapsto \sum_{j=0}^\infty e_j z^{(j)}\,, 
\] 
that satisfies for each $a \in A$ the identity 
\begin{equation} \label{Exp-Identity} 
\exp_{E,\Spm B} \circ \phi_a' \enspace = \enspace \phi_a \circ \exp_{E,\Spm B}\,. 
\end{equation} 
\end{prop} 
The proof is given in Section~\ref{SectProof}.  
As the functions $\exp_{E,\Spm B}$ are unique, they agree on overlaps, and so patching yields 
\begin{cor} 
Every abelian $A$-module $E$ on $X$ possesses a unique morphism $\exp_E: \Lie E \to E$ of rigid analytic vector bundles, such that for all $a\in A$ one has $\exp_E \circ\phi_a' =\phi_a \circ\exp_E$, and such that on every affinoid subdomain $\Spm B \subset $ of $X$ which trivializes $E$, and for all coordinates of $E|_{\Spm B}$, the morphism $\exp_E|_{\Spm B}$ takes the form given in Proposition~\ref{DefExpMap}. \end{cor} 
 
\begin{Def} 
For any {\'e}tale morphism $Y\to X$ we define 
\[ 
\Lambda(Y) \enspace := \enspace \Lambda_E(Y) \enspace := \enspace \ker \bigl(\exp_E : \enspace \Lie E(Y) \longrightarrow E(Y)\bigr) 
\] 
as the kernel of the exponential map on $Y$. This defines an {\'e}tale sheaf $\Lambda$ of $A_X$-modules on~$X$.  
\end{Def}

For an analytic point $x \in \CM(X)$ we consider the pullback $x^\ast(E,\phi)$, which is an abelian $A$-module over $\Spm k(x)$. Its exponential is obtained from $\exp_E$ by base change. The kernel of this exponential will be denoted by $\Lambda_{x^\ast E}$.


\subsection{The Key-Lemma for Abelian $A$-Modules} \label{SectFamAMod-Key} 
 
In this section we prove the abelian $A$-module analogue of our Key-Lemma~\ref{ExtendOnV}. It will allow us to show that $\Lambda$ is an overconvergent sheaf and to give in the next section a criterion for uniformizability of abelian $A$-modules.

\begin{keylemma} \label{Extendlambda} 
Let $X=\Spm B$ be affinoid and let $(E,\phi)$ be an abelian $A$-module such that $E \cong \CO_X^d$ is trivial. Let $s \in E(X)$. Let further  $x \in \CM(X)$ be an analytic point and let $\ol{\ell}$ be a given element of $(x^\ast\Lie E)\bigl(k(x)^{alg}\bigr)$, defined over an algebraic closure of $k(x)$, satisfying $\exp_E(\ol{\ell}) = s(x)$. 
 
Then there exists an {\'e}tale morphism $\pi:V\to X$ of affinoids such that $\pi V$ is a wide neighborhood of $x$ in $X$, a point $y\in \CM(V)$ above $x$, and a uniquely determined section $\ell \in \Lie E(V)$ with $\exp_E(\ell)=s|_V$ and $\ell(y) = \ol{\ell}$. Moreover, if $\ol{\ell}$ is defined over $k(x)$,  we can find $V\subset X$ as a wide affinoid neighborhood of $x$. 
\end{keylemma} 
 
\begin{proof} 
Let $\kappa: E \to \CO_X^d$ be coordinates on $E$ and let $\rho,\wt{\rho} \in L^{alg}$ with $|\rho|>|\wt{\rho}| = |\kappa'(\ol{\ell})|$. 
Consider  
\begin{eqnarray*}  
& \pi: & \DS V := \Spm B\,\langle\, \frac{z_1}{\rho}, \ldots, \frac{z_d}{\rho} \,\rangle \,/\, \bigl( \,\exp_E (z_1, \ldots, z_d) -\kappa(s)\,\bigr) \enspace \longrightarrow \enspace X \qquad \text{and}\\[0.3cm] 
&& U := \{\,y \in V: |z_i(y)| \leq |\wt{\rho}|\;,i=1, \ldots,d\,\} \enspace \subset \enspace V\,. 
\end{eqnarray*} 
Then $U$ is relatively compact in $V$ over $X$, i.e., $U \subset\subset_X V$. Furthermore, $\pi$ is {\'e}tale, since $\exp_E(z) = \sum e_j z^{(j)}$ with $e_0 = \Id$; cf.\ Proposition~\ref{DefExpMap}. By the choice of $\wt{\rho}$ the point $x$ lifts to an analytic point $y$ of $U$ with  
\[ 
\kappa'(\ol{\ell}) = \bigl(z_1(y), \ldots, z_d(y) \bigr)^T\,. 
\] 
 
We set $\ell := \kappa'{}^{-1}\,(z_1, \ldots, z_d)^T \in \Lie E(V)$ and have $\exp_E(\ell)=s|_V$ by construction. By~
\cite[Lemma 3.4.2]{JongPut} $\pi V$ is a wide neighborhood of $x$ in $X$. 
 
If $\ol{\ell}$ is defined over $k(x)$, we have $k(y) = k(x)$. By \cite[Lemma 3.1.5]{JongPut} there is an affinoid subdomain $U' \subset U$ with $y\in\CM(U')$ such that $\pi|_{U'}$ is an isomorphism $U'\to \pi(U')$. Now by \cite[Lemma 3.4.2]{JongPut} there is a wide neighborhood $V'$ of $y$ such that $\pi|_{V'}:V'\to \pi V'$ is an isomorphism and $\pi V'$ is a wide neighborhood of $x$ in $X$. 
\end{proof} 
 
\medskip 
 
The central point of the above proof may be restated in terms of Berkovich's theory of analytic spaces 
as the following result on standard {\'e}tale neighborhoods for Berkovich spaces: 
\begin{prop} 
Let $X=\Spm B$ be affinoid and let $Y:= \Spm B\langle T_1, \ldots, T_n \rangle/(f_1, \ldots, f_n)$ be {\'e}tale over $X$. Then $\CM(Y)^\circ := \{\,y\in \CM(Y): |T_i|_y <1 \enspace \text{for all} \enspace i\,\}$ is {\'e}tale over $\CM(X)$ in the sense of Berkovich. In particular, the image of $\CM(Y)^\circ$ in $\CM(X)$ is open and thus if $x$ is the image of an analytic point $y$ of $\CM(Y)^\circ$, then there is an affinoid $V$ of $Y$ whose image is a wide neighborhood of~$x$. 
\end{prop} 
 
\begin{proof} 
By Berkovich \cite[Prop.\ 3.1.4]{berkovich2} the morphism $\CM(Y)^\circ \to \CM(X)$ is quasi-finite (\cite[Def.\ 3.1.1]{berkovich2}) and hence {\'e}tale (\cite[Def.\ 3.3.4]{berkovich2}). The last statement follows from \cite[Prop.\ 3.2.7]{berkovich2}. 
\end{proof} 
 
\begin{remark} 
The above result is optimal in that one cannot expect $\CM(Y)$ to be Berkovich-{\'e}tale over $\CM(X)$. For example consider $X=\Spm L\langle T\rangle$ and $Y= \Spm L\langle T/\pi\rangle$ for $|\pi|<1$. The morphism $\CM(Y) \to \CM(X)$ is not quasi-finite on the boundary of~$\CM(Y)$.  
\end{remark} 
  
\medskip 
 
After this digression to Berkovich spaces, let us return to the theory of abelian $A$-modules. As a first application of the above we obtain the following analogue to Corollary~\ref{SatUModTau}. 
 
\begin{cor} \label{SatUModTMod} 
Let $(E,\phi)$ be an abelian $A$-module on a connected $X$. 
If $\pi:Y \to X$ is a general morphism (cf.~Definition~\ref{DefGenMorph}), then the map 
\[ 
\Lambda_E(X) \to \Lambda_{\pi^\ast \!E}(Y)\,,\qquad \lambda \mapsto \pi^\ast \lambda 
\] 
is injective and its image is a saturated $A$-submodule in~$\Lambda_{\pi^\ast E}(Y)$. 
\end{cor} 
 
\begin{proof} 
If $X$ and $Y$ are spectra of algebraically closed complete fields $K$ and $L$, the assertion is clear since then we have $\Lambda_E(K)\cong \Lambda_{\pi^\ast E}(L)$. As in the proof of Corollary~\ref{SatUModTau}, it therefore suffices to show the assertion for arbitrary $X$ and for $Y=\Spm \wh{k(x)^{alg}}$ for any $x\in\CM(X)$, which from now on we assume. 
 
Since $\lambda\in \Lambda_E(X)$ satisfies the {\'e}tale equation $\exp_E(\lambda)=0$ and $X$ is connected, the morphism 
\[ 
h: \Lambda_E(X) \longrightarrow \Lambda_{x^\ast E}\bigl(\wh{k(x)^{alg}}\bigr)\,,\quad \lambda \mapsto \lambda(x) 
\] 
is injective. To prove that the image is saturated, let $\lambda \in \Lambda_E(X)$, $\ol{\mu} \in \Lambda_{x^\ast E}\bigl(k(x)^{alg}\bigr)$ and $a\in A$, $a \neq 0$, satisfying $\lambda(x)=\phi'_a(\ol{\mu})$. Since $\phi'_a$ is an automorphism of $\Lie E$, $\ol{\mu}$ is defined over $k(x)$. We let $\mu := \phi'_a{}^{-1}(\lambda) \in \Lie E(X)$ and claim that $\mu\in \Lambda_E(X)$. 
 
Choosing an admissible affinoid covering $\Spm B_i$ of $X$ we can trivialize $E$. There is an $i$ such that $x$ lies in $\Spm B_i$. By Lemma~\ref{Extendlambda} ($s=0$) there is a wide affinoid neighborhood $V\subset \Spm B_i$ of $x$ and a section $\wt{\mu} \in \Lambda_E(V)$ with $\wt{\mu}(x) = \ol{\mu}$. By the injectivity of $h$ we have $\lambda|_V = \phi'_a(\wt{\mu})$, i.e., $\mu|_V=\wt{\mu}$. This implies, that $\exp_E(\mu)\,|_V = 0$. Since $E$ is coherent and $X$ is connected, $E(X)\to E(V)$ is injective, and we conclude that $\exp_E(\mu) = 0$ on $X$ as required. 
\end{proof} 
 
\medskip 
 
\begin{cor} 
Let $(E,\phi)$ be an abelian $A$-module on a connected $X$. Then 
\begin{enumerate} 
\item 
For every analytic point $x\in\CM(X)$ the {\'e}tale sheaves $x^\ast\Lambda_E$ and $\Lambda_{x^\ast E}$ on $\Spm k(x)$ are canonically isomorphic. 
\item  
For $y\in \CM_\et(X)$ the stalk of $\Lambda_E$ at $y$ is given by 
  $\Lambda_{y^\ast E}\bigl(k(y)\bigr)$. 
\item  
$\Lambda_E$ is overconvergent on the small {\'e}tale site of~$X$. 
\item  
The homomorphism $\Lie E\to E$ is surjective if and only if it is so  
  for the fibers $y^\ast E$ at all {\'e}tale analytic points $y\in \CM_\et(X)$. 
\end{enumerate} 
\end{cor} 
The proof is analogous to that of Corollary~\ref{StalkOfFtauAndOverC}, where instead of  
Lemma~\ref{ExtendOnV} and Corollary~\ref{SatUModTau},  
one uses Lemma~\ref{Extendlambda} and Corollary~\ref{SatUModTMod}, and so it is omitted. 
 


\subsection{Uniformizability of Abelian $A$-Modules} \label{SectFamAMod-Unif}

We now generalize Anderson's notion and characterization of 
uniformizability of an abelian $A$-module on a point 
\cite[Thm. 4]{anderson} to the case of an abelian $A$-modules on a 
general base~$X$. 
 
\begin{thm} \label{UnivAMod} 
Let $E$ be an abelian $A$-module of dimension $d$ and rank $r$ on $X$. Then the following are equivalent. 
\begin{enumerate} 
\item 
The rigid analytic $\tau$-sheaf $\UCF$ over $A(1)$ on $X$ associated to $E$ is analytically trivial, 
\item 
there exists a connected temperate {\'e}tale Galois covering  $Z$ of $X$ such that $\Lambda(Z)$ is a projective $A$-module of rank $r$. 
\item 
there exists a covering for the {\'e}tale topology $\{\,Y_i \to X\,\}$ such that $\Lambda(Y_i)$ is a projective $\ul{A}_X(Y_i)$-module of rank $r$ for all $i$, 
\item 
for every analytic point $x \in \CM(X)$ the abelian $A$-module $x^\ast E$ is uniformizable in the sense of Anderson, 
\item 
the following short sequence of {\'e}tale sheaves on $X$ is exact 
\[ 
0 \longrightarrow \Lambda \longrightarrow \Lie E \xrightarrow{\exp_E} E \longrightarrow 0\,. 
\] 
\end{enumerate} 
\end{thm}  
 
If the equivalent conditions of Theorem~\ref{UnivAMod} are satisfied we call $E$ {\em uniformizable}. 

\begin{proof} 
$(a) \Rightarrow (b)$. By Theorem~\ref{PointwiseTau} there exists a connected temperate {\'e}tale Galois covering $\pi:Z\to X$, such that $\UCF^\tau(Z)$ is a projective $A$-module of rank $r$. 
As in Anderson one shows that the action $A \to \End(\Lie E)$ extends uniquely to a continuous action of $K_\infty$ and that the pairing of $A$-modules 
\[ 
\Lambda(Z) \times \UCF^\tau(Z) \longrightarrow \Omega_A\,,\quad (\lambda,f) \mapsto \omega_{\lambda, f}\,, 
\] 
which is defined by requiring that $f(\exp_E(z \cdot\lambda)) = \Res_\infty(z\,\omega_{\lambda,f})$ for all $z\in K_\infty$, is perfect. For more details see \cite[{\S} 9.2]{boeckle}. We conclude that $\Lambda(Z)$ is a projective $A$-module of rank~$r$.

The implications $(b) \Rightarrow (c) \Rightarrow (d) \Rightarrow (a)$ are clear by Corollary~\ref{SatUModTMod}, Anderson \cite[Thm. 4]{anderson}, and Theorem~\ref{PointwiseTau}.

$(d) \Rightarrow (e)$. We have to show that $\exp_E$ is surjective, i.e., that for every {\'e}tale morphism of rigid analytic spaces $Y\to X$ and every $s \in E(Y)$ there exists a covering for the {\'e}tale topology $\{V_i \to Y\}_{i\in I}$ and for each $i$ an element $\ell_i \in \Lie E(V_i)$ satisfying $\exp_E(\ell_i) = s|_{V_i}$. By passing to an admissible covering we may assume that $Y=\Spm B$ is a connected affinoid and $E|_Y$ is trivial. By assumption (d) and Lemma~\ref{Extendlambda} every analytic point $y\in \CM(Y)$ admits a wide affinoid neighborhood $U_y$ and an {\'e}tale surjective morphism $V_y \to U_y$, such that there exists an $\ell_y \in \Lie E(V_y)$ with $\exp_E(\ell_y) = s|_{V_y}$. Since $\CM(Y)$ is compact, finitely many of the $U_y$ suffice, giving  the desired covering.

$(e) \Rightarrow (d)$. Let $x \in \CM(X)$ be an analytic point. First let $\ol{s} \in x^\ast E \bigl(k(x)^{sep}\bigr)$ be defined over a separable closure $k(x)^{sep}$ of $k(x)$. As $E$ is coherent the section $\ol{s}$ is actually defined over a finite separable extension $L'$ of $k(x)$. By~
\cite[Remark  2.1.2]{JongPut} there is an {\'e}tale morphism $Y \to X$ of affinoids and an analytic point $y$ of $Y$ above $x$ with $k(y)=L'$. We extend $\ol{s}$ to a section $s\in E(Y)$. Then our assumption (e) implies that there is an {\'e}tale morphism $V\to Y$ whose image is a special subset $W$ with $y\in\CM(W)$ and an $\ell\in \Lie E(V)$ with $\exp_E(\ell) = s|_V$. Therefore $\ol{\ell} :=\ell(x)$ satisfies $\exp_E(\ol{\ell}) = \ol{s}$. 
 
Since $k(x)^{sep}$ is dense in $k(x)^{alg}$ by \cite[3.4.1/6]{BGR} and $\exp_E$ is a local isomorphism we conclude that $\exp_E :(x^\ast\Lie E)\bigl(k(x)^{alg}\bigr) \to (x^\ast E)\bigl(k(x)^{alg}\bigr)$ is surjective as required. 
\end{proof}

 
\subsection{The Proof of Proposition~\ref{DefExpMap}}\label{SectProof} 
 
\begin{proof} 
Denote by $B\{\{\sigma\}\}$ the ring of formal series $\sum_{k=0}^\infty b_k\sigma^k$, $b_k\in B$, under component wise addition and composition as multiplication. With this definition, we have $\phi_a\in M_d(B\{\{\sigma\}\})$. In a first step we construct a formal function $\exp_{E,\Spm B}\in M_d(B\{\{\sigma\}\})$ satisfying (\ref{Exp-Identity}). For the argument, we follow Goss \cite[p.~75f]{Goss} which ultimately goes back to Drinfeld. Basically the argument goes by first constructing $\exp_{E,\Spm B}$ which satisfies (\ref{Exp-Identity}) for a single $a\in A\smallsetminus \BF_q^{\,alg}$ and then one shows that the same formal function works for any $a\in A$. 
 
Let $a\in A$ be a non-constant element. Because $G_0(a)-a$ is nilpotent, and $a\in A\subset K_\infty\subset L$ is a unit in $L$, we have $G_0(a)\in \GL_d(B)$. A recursive argument, based on the formula $(1+x)^{-1}=\sum x^n$, shows that in fact $\phi_a$ is invertible in $M_d(B\{\{\sigma\}\})$. If follows that $\phi$ may be regarded as a homomorphism  
$$\phi\!:A\to M_d(B\{\{\sigma\}\}),$$ 
with $\phi(A\smallsetminus \{0\})\subset \GL_d(B\{\{\sigma\}\})$. Let us formulate and prove the following lemma:

\begin{lemma} 
Let $f=\sum a_k\sigma^k\in M_d(B\{\{\sigma\}\})$. Assume that $a_0\in M_d(B)$ can be written as $\alpha \Id+N_1$ for some $\alpha\in L\smallsetminus \BF_q^{\,alg}$ and some nilpotent matrix $N_1$ and suppose that we have $e_0\in M_d(B)$ with $N_1e_0=e_0N_2$ for some nilpotent matrix $N_2$. Then there exists a unique $\lambda_f=\sum e_j\sigma^j\in M_d(B\{\{\sigma\}\})$ such that  
\begin{equation}\label{LambdaFformula} \lambda_f(\alpha \Id+N_2)=f \lambda_f.\end{equation} 
\end{lemma} 
\begin{proof} 
The conditions on the coefficients $e_j$ of $\lambda_f$ in formula (\ref{LambdaFformula}) can, by sorting the expression according to powers of $\sigma$, be expressed equivalently in the following recursive way: 
\begin{equation}\label{ConditionEj} 
\forall j\in\BN \,:\,e_j(\alpha^{(j)}+N_2^{(j)})=(\alpha+N_1)e_j +\sum_{k=1}^{j} a_k e_{j-k}^{(k)}. 
\end{equation} 
Because $N_1$ and $N_2$ are nilpotent and $\alpha^{(j)}-\alpha$ is a unit $(\alpha\notin \BF_q^{\,alg})$ for $j>0$, the linear operator  
$$T_j\!:e\mapsto e(\alpha^{(j)}+N_2^{(j)})-(\alpha+N_1)e$$ 
on matrices $e\in M_d(B)$ is invertible. By recursively solving for the $e_j$ the lemma follows. 
\end{proof} 
We now apply the lemma to $f=\phi_a$ and $N_1=N_2$, and obtain a unique $\lambda_a\in M_d(B\{\{\sigma\}\})$ whose constant term is the identity matrix $\Id$ and such that 
$\lambda_a G_0(a)=\phi_a \lambda_a$. Since the constant term is invertible, $\lambda_a$ is a unit in $M_d(B\{\{\sigma\}\})$, and so we find 
$$\phi_a=\lambda_a G_0(a) \lambda_a^{-1}.$$ 
As $\phi(A)\subset M_d(B\{\{\sigma\}\})$ is commutative, the following is an important observation: 
\begin{lemma} 
Let $f$ be as in the previous lemma, and assume that $N:=N_1=N_2$ and $e_0\in\GL_d(B)$. Then for the centralizer of $f$ in $M_d(B\{\{\sigma\}\})$ we have 
$$\Cent_{M_d(B\{\{\sigma\}\})}(f)=\lambda_f\Cent_{M_d(B)}(\alpha \Id+N)\lambda_f^{-1}.$$ 
\end{lemma} 
\begin{proof} 
As $f=\lambda_f(\alpha \Id+N)\lambda_f^{-1}$, any element on the right clearly commutes with $f$. Therefore we only have to show the inclusion~$\subseteq$. Suppose that we are given $g\in \Cent_{M_d(B\{\{\sigma\}\})}(f)$. It is then clear that $\tilde g:=\lambda_f^{-1}g\lambda_f$ lies in $\Cent_{M_d(B\{\{\sigma\}\})}(\alpha \Id+N)$. 
If we write $\tilde g=\sum d_j\sigma^j$, then this is equivalent to $$\forall j>0\,:\, T_j d_j =d_j(\alpha^{(j)}+N^{(j)})-(\alpha+N)d_j=0$$and $d_0(\alpha \Id+N)=(\alpha \Id+N)d_0$. Since the operators $T_j$ from the proof of the previous lemma are linear and invertible, we conclude $d_j=0$ for $j>0$ and $d_0\in \Cent_{M_d(B)}(\alpha \Id+N)$, as asserted. 
\end{proof}

By the lemma and the commutativity of $\phi(A)$, we have $\lambda_a^{-1}\phi_{a'}\lambda_a\in \Cent_{M_d(B)}\bigl(G_0(a)\bigr)$ for any $a'\in A$. Because $\lambda_a$ has constant term $\Id$ and $\phi_{a'}$ has constant term $G_0(a')$ we have in fact 
$$\lambda_a^{-1}\phi_{a'}\lambda_a=G_0(a')$$ 
for all $a'\in A$. Thus if we fix an element $t\in A\smallsetminus \BF_q^{\,alg}$ and define $\exp_{E,\Spm B}:=\lambda_t = \sum_j e_j \sigma^j$, then $\exp_{E,\Spm B}$ is the unique element in $M_d(B\{\{\sigma\}\})$ with constant term $\Id$ which satisfies the condition 
$$ \forall a\in A:\quad \exp_{E,\Spm B}G_0(a)\es=\es\phi_a\exp_{E,\Spm B}\es\in M_d(B\{\{\sigma\}\}).$$ 
 
\medskip 
 
The proof of the analyticity of $\exp_{E,\Spm B}$ proceeds as in Anderson \cite[{\S} 2.1]{anderson}, except that we can not assume that $G_0(t)$ is in Jordan canonical form. So we need finer estimates to ensure the convergence of $\exp_{E,\Spm B}$. We fix a complete $L$-algebra norm $|\;\;|$ on the affinoid algebra $B$. For any matrix $\Delta$ with entries in $B$ we let $|\Delta|$ be the maximum over the norms of the entries. We denote by $\theta$ the image of $t$ in $B$. Observe that $|\theta|>1$. The matrix $N:=G_0(t)-\theta\Id$ is nilpotent, say with $N^{m+1}=0$. But we don't have $|N|\leq 1$. We set $G_k=G_k(t)$. The identity (\ref{ConditionEj}) for $a=t$ implies that the $e_j$ are uniquely determined by the recursion  
\[ 
e_j \enspace = \enspace \sum_{k=1}^s \,\sum_{i=0}^{2m}\,\sum_{\nu=0}^i \,\frac{\binom{i}{\nu}}{\bigl(\theta^{q^j}-\theta\bigr)^{i+1}} \;N^\nu \;G_k \;e_{j-k}^{(k)} \;\bigl(N^{(j)}\bigr)^{i-\nu}\,. 
\] 
Thus by induction we see that 
\[ 
e_j \enspace = \enspace \sum_{\ul{k}} \,\sum_{\ul{i}} \,\sum_{\ul{\nu}} \,\frac{\binom{\ul{i}}{\ul{\nu}}}{\theta^{^{\SC q^j(|i|+l)}}(1+\epsilon)} \;H_j^{(0)}\cdot \ldots \cdot H_1^{(j-1)} \;e_0^{(j)} \;\bigl(N^{(j)}\bigr)^{|i|-|\nu|}\,, 
\] 
where the first, second and third sum runs over  
\begin{eqnarray*} 
&&\{\ul{k} = (k_1, \ldots, k_j) \in \Nz_0^j: k_\mu \leq s,\; k_j\neq 0,\; \text{ and } k_1+ \ldots+k_\mu = \mu \text{ if } k_\mu \neq 0\}\,,\\[0.2cm] 
&&\{\ul{i} = (i_1, \ldots, i_j) \in \Nz_0^j: \enspace i_\mu \leq 2m, \;\text{ and }i_\mu = 0 \text{ if } k_\mu=0\}\,,\\[0.2cm] 
&&\{\ul{\nu} = (\nu_1, \ldots, \nu_j) \in \Nz_0^j: \enspace \nu_\mu \leq i_\nu\}\,, 
\end{eqnarray*} 
respectively, where we abbreviate  
\begin{eqnarray*} 
&&\binom{\ul{i}}{\ul{\nu}} \enspace = \enspace \binom{i_1}{\nu_1}\cdot \ldots \cdot\binom{i_j}{\nu_j}\,,\\[0.2cm] 
&&|i| \enspace = \enspace i_1+\ldots +i_j\,,\qquad |\nu|\enspace = \enspace \nu_1+\ldots+\nu_j\,,\\[0.2cm] 
&&l \enspace = \enspace l(\ul{k})\enspace = \enspace \#\{\mu: \;k_\mu \neq 0\}\,,\\[0.2cm] 
&&H_\mu \enspace = \enspace H_\mu(\ul{k},\ul{\nu}) \enspace = \enspace \left\{ 
\begin{array}{l@{\quad\text{if}\quad}l} 
\Id & k_\mu = 0\,,\\ 
N^{\nu_\mu} G_{k_\mu} & k_\mu\neq 0\,, 
\end{array}\right. 
\end{eqnarray*} 
and where $\epsilon=\epsilon(\ul{k},\ul{i},\ul{\nu})$ is a polynomial in $\theta^{-1}$ with $|\epsilon|<1$. Now the constraints on $\ul{k}$ imply that $l \geq j/s$. Hence the absolute value of the fraction is at most 
\[ 
|\theta|^{-q^j l} \leq |\theta|^{-q^j j/s}\,. 
\] 
Furthermore, if $|i|-|\nu|>m$, the most right factor is zero. Let $M$ be a common upper bound for $|G_\mu|$, $|N|, \ldots,|N^m|$. Thus we see that 
\[ 
|H_\mu^{(j-\mu)}| \leq M^{2q^{j-\mu}} \qquad \text{and}\qquad |\bigl(N^{(j)}\bigr)^{|i|-|\nu|}| \leq M^{q^j}\,. 
\] 
Thus we conclude as required 
\[ 
q^{-j} \log |e_j| \enspace \leq \enspace -j/s\,\log|\theta| + q^{-j}\frac{q^j-1}{q-1} \log M^2 + \log M \enspace \longrightarrow \enspace -\infty \quad (j\to \infty)\,. 
\] 
This proves Proposition~\ref{DefExpMap}. 
\end{proof}

\begin{appendix}

 
\section{Review of Analytic Spaces and {\'E}tale Sheaves}  
\setcounter{equation}{0}\label{Review}


\subsection{Analytic Points}\label{AnalyticPoints}

Given a category of sheaves of abelian groups on a space $X$ one wants to dispose of a sufficiently large set of points to detect exactness via stalks at these points. On a rigid analytic space $X$ the {\em classical points\/} (the ones usually considered in rigid analytic geometry, whose residue fields are finite extensions of the base field) are too few if one works with sheaves which are not coherent. To remedy this, van der Put \cite{Put} introduced the larger class of {\em analytic points\/} together with {\em overconvergent sheaves\/} (Definitions~\ref{AnalyticPoint} and \ref{DefOverconvergent} below). In the category of overconvergent sheaves on $X$ a sequence is exact if and only if for every analytic point $x$ of $X$ the associated sequence of stalks at $x$ is exact; cf.~
\cite[{\S} 2]{schneider} or~
\cite[{\S} 2]{JongPut}. 
 
An important tool in the study of uniformizability of $\tau$-sheaves are their associated sheaves of $\tau$-invariants. The latter will turn out to be overconvergent ({\'e}tale) sheaves; cf.\ Corollary~\ref{StalkOfFtauAndOverC}. We therefore now briefly recall the concepts of analytic point, overconvergent ({\'e}tale) sheaf, etc., and describe some basic properties. For more details see the articles \cite{berkovich}, \cite{JongPut}, \cite{Put}, \cite{PutSchneider}, \cite{schneider} by Berkovich, de Jong, van der Put, and Schneider.

\bigskip

Let $L$ be a complete, non-archimedean valued field and denote by $X = \Spm B$ an affinoid rigid analytic space over $L$. Following \cite{JongPut}, we call an admissible subset $V\subset X$ {\em special\/} if $V=\cup V_i$ is a finite union of affinoid subdomains $V_i\subset X$.

\begin{Def} \label{AnalyticPoint} 
An {\em analytic point\/} $x$ of $X$ is a semi-norm $|\,.\,|_x:B \to \Rz_{\geq 0}$ which satisfies: 
\begin{enumerate} 
\item 
$|f+g|_x \leq \max\{\,|f|_x,|g|_x\,\}$ \quad for all $f,g \in B$, 
\item  
$|fg|_x = |f|_x \,|g|_x$ \quad for all $f,g \in B$, 
\item 
$|\lambda|_x = |\lambda|$ \quad for all $\lambda \in L$, 
\item 
$|\,.\,|_x:B \to \Rz_{\geq 0}$ is continuous with respect to the norm topology on $B$. 
\end{enumerate} 
The set of all analytic points of $X$ is denoted~$\CM(X)$. 
\end{Def}

On $\CM(X)$ one considers the coarsest topology such that for every $f\in B$ the map $\CM(X) \to \Rz$ given by $x \mapsto |f|_x$ is continuous. Equipped with this topology, $\CM(X)$ is a compact Hausdorff space; cf.~
\cite[Thm.\ 1.2.1]{berkovich}, 
\cite[{\S} 1]{Put} or~
\cite[{\S} 1]{schneider}.  
 
Every classical point $x \in X$ defines a semi-norm as above via 
\[ 
B \longrightarrow k(x) \xrightarrow{\enspace |\,.\,|_x\;} \Rz_{\geq 0} 
\] 
where $|\,.\,|_x$ is the unique norm on its residue field $k(x)$ which extends the norm of $L$. Thereby one obtains an injective map $X \to \CM(X)$. With respect to the canonical topology on $X$, cf.~
\cite[{\S} 7.2.1]{BGR}, this map is a homeomorphism onto a dense subset in $\CM(X)$; cf.~
\cite[Prop.\ 2.1.15]{berkovich}.  
 
For a classical point $x\in X$ the residue field $k(x)=B/\ker|\,.\,|_x$ is a finite extension of $L$ and thus complete with respect to $|\,.\,|_x$. 
For an arbitrary analytic point $x \in \CM(X)$ this is no longer true. In this situation one defines $k(x)$ to be the completion of the fraction field of $B/\ker |\,.\,|_x$ with respect to $|\,.\,|_x$ and calls it the {\em (complete) residue field of $x$}. This gives rise to a continuous homomorphism $B \to k(x)$ of $L$-algebras. Hence the analytic points coincide with the equivalence classes of analytic points defined by Schneider \cite[{\S} 1]{schneider}.  
 
Every morphism $f:Y\to X$ of affinoid rigid analytic spaces over $L$ induces a continuous morphism $\CM(f):\CM(Y) \to \CM(X)$ by mapping the semi-norm $\CO(Y) \to \Rz_{\geq 0}$ to the composition $\CO(X) \to \CO(Y) \to \Rz_{\geq 0}$. In particular, for an affinoid subdomain $U\subset X$ this morphism identifies $\CM(U)$ with a closed subset of $\CM(X)$. 
 
\begin{Def} \label{DefM(U)} 
If $U\subset X$ is an admissible subset with an admissible affinoid covering $\{U_i\}$ one defines 
\[ 
\CM(U) \enspace := \enspace \bigcup_i \CM(U_i) \enspace \subset \enspace \CM(X)\,. 
\] 
\end{Def} 
 
It follows from 
\cite[Lemma 1.3]{schneider} that this definition is independent of the chosen covering $\{U_i\}$. In particular, if $U$ is itself affinoid one recovers Definition~\ref{AnalyticPoint}.

\begin{Def} \label{DefWideNbhd} 
An admissible subset $U\subset X$ is called a {\em wide neighborhood\/} of an analytic point $x\in \CM(X)$ if $x$ lies in the interior $\CM(U)\open$. 
\end{Def} 
 
Note that if $x\in X$ is a classical point then every admissible $U\subset X$ containing $x$ is a wide neighborhood of $x$. However, for general analytic points this is not the case.

\begin{Def} \label{DefRelCompact} 
For admissible subsets $V\subset U\subset X$ one calls {\em $U$ a wide neighborhood of $V$ in $X$} and writes $V\subset\subset_X U$ if $\CM(V)$ is contained in the interior $\CM(U)\open$. 
\end{Def} 
 
If $V$ and $U$ are affinoid, Kiehl~\cite[Def.\ 2.1]{kiehlEndl} originally used the symbol $V\subset \subset_X U$ to mean that $V$ is {\em relatively compact in $U$ over $X$\/}, i.e. that there is an affinoid generating system $f_1,\ldots,f_r$ of $\CO_X(U)$ over $\CO_X(X)$ such that 
\[ 
V\enspace \subset \enspace \{x\in U: \enspace|f_i(x)|<1 \quad\text{for }1\leq i\leq r\}\,. 
\] 
By the following result, Kiehl's definition is equivalent to the one given in Definition~\ref{DefRelCompact}.

\begin{prop} 
Let $V\subset U\subset X$ be admissible subsets. 
\begin{aufz} 
\item 
If $V$ and $ U$ are affinoid then $U$ is a wide neighborhood of $V$ in $X$ if and only if $V$ is relatively compact in $U$ over $X$. 
\item 
If $V=\cup V_i$ is special then $V\subset\subset_X U$ if and only if $V_i\subset\subset_X U$ for every $i$. 
\item 
If $V$ is special then $V\subset\subset_X U$ if and only if there exists a special subset $W\subset U$ such that $V\subset\subset_X W$. 
\item 
$\DS\CM(U)\open \enspace = \enspace \bigcup \{\CM(V): V \text{ affinoid with } V\subset\subset_X U\}$\,. 
\end{aufz} 
\end{prop} 
 
\begin{proof} 
Part (a) is from 
~\cite[Prop.\ 1.23]{schneider}, part~(b) follows immediately from the definition, and (c) and (d) are 
from
~\cite[Remark 1.24]{schneider} and
~\cite[Lemma 1.12]{schneider}, respectively. 
\end{proof} 
 
The following two propositions are taken from \cite[Prop.\ 3.3 and Lemma 3.1]{schneider}, and \cite[Prop.\ 3.5]{schneider}, respectively. 
 
\begin{prop} 
For a special subset $V\subset X$, the set $\CM(X\smallsetminus V) = \CM(X) \smallsetminus \CM(V)$ is open in~$\CM(X)$. 
\end{prop} 
 
\smallskip

\begin{prop} \label{AdmCovering} 
Let $\DS U = \bigcup\{U_i:i\in I\}$ be a covering of an admissible subset $U \subset X$ by admissible subsets $U_i\subset X$ with $\CM(U_i) \subset \CM(X)$ open for all $i$. Then this covering is admissible if and only if $\DS \CM(U) = \bigcup\{\CM(U_i):i\in I\}$. In this case $\CM(U)$ is open in $\CM(X)$. 
\end{prop}

\medskip 
 
We further need the following lemma. 
 
\begin{lemma} \label{IntersecOpen} 
Let $X$ be an affinoid rigid analytic space and let $X_i,\,i\in I$ be a finite affinoid covering of $X$. Let $M \subset \CM(X)$ be a subset. Then $M$ is open in $\CM(X)$ if and only if $M \cap \CM(X_i)$ is open in $\CM(X_i)$ for every $i$. 
\end{lemma}

\begin{proof} 
This is a standard argument in point set topology, which we repeat for sake of completeness: Define $M^c:=\CM(X)\smallsetminus M$. Then we need to show that $M^c$ is closed in $\CM(X)$ if and only if $M^c\cap \CM(X_i)$ is closed for every $i$. This however is clear since the compact sets $\CM(X_i)$ form a finite cover of the Hausdorff space $\CM(X)$ by closed subsets. 
\end{proof}

\medskip

The notion of analytic point can be defined for general rigid analytic spaces $X$ over $L$. Namely, if $\{\,X_i\,\}_{i\in I}$ is an admissible affinoid covering of $X$, the analytic points of $X$ are just the analytic points of the $X_i$ modulo the obvious identification coming from the inclusions $X_i\cap X_j \subset X_i,X_j$. See~
\cite[{\S} 2]{schneider} for a precise  definition. The definition of a suitable topology on $\CM(X)$, so that the resulting space is Hausdorff again, is more subtle.

If $X$ is quasi-separated and admits a locally finite admissible affinoid covering $\{X_i\}$ this can be done as follows: (Recall that $X$ is called {\em quasi-separated\/} if the intersection of any two admissible affinoid subsets of $X$ is a finite union of affinoid subsets and that a covering $\{X_i\}$ is called {\em locally finite} if each $X_i$ meets only finitely many $X_j$.) One defines 
\[ 
\CM(X) \enspace := \enspace \bigcup \CM(X_i) \;/\sim 
\] 
where $\sim$ is the equivalence relation obtained by gluing $\CM(X_i)$ and $\CM(X_j)$ along the subset $\CM(X_i\cap X_j)$. The topology on $\CM(X)$ is the finest for which all the natural inclusions $\CM(X_i)\into\CM(X)$ are continuous. For such $X$ Lemma~\ref{IntersecOpen} is still true; cf.~
\cite[Lemma 5.3]{PutSchneider}.

For the general situation V.G. Berkovich in \cite{berkovich}, \cite{berkovich2}, has developed a 
theory of analytic spaces using the compact spaces $\CM(X)$ for 
affinoid $X$ as building blocks. However, the category of analytic 
spaces he obtains is different from the category of all rigid analytic 
spaces. As we want to remain within the framework of rigid analytic 
spaces, we therefore speak of $\CM(X)$ only for rigid analytic spaces 
$X$ which are quasi-separated and paracompact. Recall that a 
rigid analytic space is called {\em paracompact\/} if every admissible 
covering admits a locally finite refinement.

 
\subsection{General Morphisms} \label{GenMorph} 
 
We recall the notion of {\em general morphisms\/} introduced by de Jong, van der Put \cite[{\S} 2.6]{JongPut}. 
Consider an extension of complete valued fields $L \subset L'$. In~
\cite[{\S} 9.3.6]{BGR} there is constructed a base change functor $X \mapsto X \widehat{\otimes}_L L'$ from quasi-separated rigid analytic $L$-spaces to quasi-separated rigid analytic $L'$-spaces.

\begin{Def}\label{DefGenMorph} 
Let $L \subset L'$ be an extension of complete valued fields and let $X$ and $Y$ be rigid analytic spaces over $L$ and  $L'$, respectively, with $X$ quasi-separated. A {\em general morphism\/} $f:Y \to X$ is a morphism $f: Y \to X \widehat{\otimes}_L L'$ of rigid analytic spaces over $L'$. 
\end{Def}

If both $X$ and $Y$ are affinoid, then this is simply a continuous homomorphism $\CO(X) \to \CO(Y)$, since any such factors as  $\CO(X) \to \CO(X)\widehat{\otimes}_L L' \to \CO(Y)$. If $Y\to X$ is a general morphism, and $Z\to X$ is a morphism of rigid analytic spaces over $L$, then we can form the fiber product 
\[ 
Y \times_X Z \enspace := \enspace Y \times_{X \widehat{\otimes}_L L'} Z \widehat{\otimes}_L L'\,. 
\] 
Every general morphism $f:Y \to X$ gives rise to a pullback functor $f^\ast$ from (quasi-)coherent $\CO_X$-modules to (quasi-)coherent $\CO_Y$-modules. 
Any analytic point $x \in \CM(X)$ can be viewed as a general morphism $i:\Spm k(x) \to X$.

 
\subsection{{\'E}tale Sheaves on Rigid Analytic Spaces}\label{EtaleSheaves}

A morphism $f:Y\to X$ of rigid analytic spaces over $L$ is called {\em {\'e}tale} if for every (classical) point $y\in Y$ the induced homomorphism of local rings $\CO_{X,f(y)} \to \CO_{Y,y}$ is flat and unramified. See~
\cite[{\S} 3]{JongPut} for a thorough discussion of this notion. 
 
Let $X$ be a rigid analytic space over $L$. We recall the definition of the {\'e}tale site of $X$ from Schneider, Stuhler \cite[p. 58]{SchnStuhler}; cf.\ also~
\cite[{\S} 3.2]{JongPut}. The underlying category of the site $X_{\acute{e}tale}$ is the category of all {\'e}tale morphisms $f: Y\to X$ of rigid analytic spaces over $L$. A morphism from $f$ to $f'$ is a morphism $g:Y \to Y'$ such that $f'\circ g=f$. The morphism $g$ is automatically {\'e}tale. 
 
\begin{Def} \label{DefEtCov} 
A family of {\'e}tale morphisms $\{\,g_i:Z_i \to Y\,\}_{i\in I}$ is a {\em covering for the {\'e}tale topology\/} if it has the following property: 
 
\medskip 
 
\noindent 
\hspace{0.05\textwidth} 
\parbox{0.9\textwidth}{ 
For every (some) choice of admissible affinoid covering $Z_i = \bigcup_j Z_{i,j}$ one has $Y = \bigcup_{i,j} g_i(Z_{i,j})$, and this is an admissible covering in the Grothendieck-topology of $Y$.}  
\end{Def} 
 
\medskip 
 
Clearly any admissible covering of $Y$ is a covering for the {\'e}tale topology. 
 
The property in Definition~\ref{DefEtCov} is local on $Y$ in the following sense: if $Y=\bigcup Y_l$ is an admissible affinoid covering, then $\{\,g_i:Z_i \to Y\,\}$ is a covering for the {\'e}tale topology if and only if for all $l$ the same is true for the covering $\{\,g_i:g_i^{-1}(Y_l) \to Y_l\,\}$. This implies that if $\{\,Z_i \to Y\,\}$ and $\{\,W_{i,j} \to Z_i\,\}$ for all $i$ are coverings for the {\'e}tale topology, then $\{\,W_{i,j} \to Y\,\}$ is a covering for the {\'e}tale topology.  
 
\smallskip 
 
The category $X_{\acute{e}tale}$ equipped with the family of  coverings for the {\'e}tale topology is thus a site, called the {\em {\'e}tale site of $X$}. The sheaves on this site are called {\em {\'e}tale sheaves on $X$}. 
 
If $\{\,Z_i \to Y\,\}$ is a covering for the {\'e}tale topology, and $Y' \to Y$ is a general morphism, then the fiber product $\{\,Y' \times_Y Z_i \to Y'\,\}$ is a covering for the {\'e}tale topology, cf.~
\cite[Lemma 3.2.1]{JongPut}). So every general morphism $f:X' \to X$ induces a morphism of sites $X'_{\acute{e}tale} \to X_{\acute{e}tale}$. 
 
\begin{examples}\label{ExsOfEtSh} 
The following are examples of {\'e}tale sheaves on $X$ which we will need in the sequel ($f:Y\to X$ will denote a general object of $X_{\acute{e}tale}$): 
\begin{aufz} 
\item 
The structure sheaf $\Gr_a$ defined by $Y \mapsto \Gamma(Y,\CO_Y)$. 
\item 
For any quasi-coherent sheaf $\CF$ of $\CO_X$-modules on $X$ we define the the {\'e}tale sheaf $W(\CF)$ on $X_{\acute{e}tale}$ by $Y\mapsto \Gamma(Y,f^\ast\CF)$, where $f^\ast$ denotes the pullback of quasi-coherent modules. In particular, one has $W(\CO_X)\cong\BG_a$. Any {\'e}tale sheaf $W(\CF)$ is a sheaf of $\Gr_a$-modules. 
\item 
Any representable sheaf $Y \mapsto \Hom_X(Y,Z)$ given by some rigid analytic space $Z$ over $X$. 
\item 
For any group or ring $B$ the constant {\'e}tale sheaf $\nul{B}_X$ is defined by $Y \mapsto \prod_{\pi_0(Y)}B$, where $\pi_0(Y)$ is the set of connected components of $Y$. (The restriction maps are the obvious ones.) The sheaf $\nul{B}_X$ is in fact representable, namely by $\coprod_{b\in B}X$. If $X$ is clear from the context we will also write $\nul{B}$ instead of $\nul{B}_X$. 
\end{aufz} 
\end{examples}

\begin{Def}\label{DefOverconvergent} 
A sheaf $\CS$ on $X_{\acute{e}tale}$ is called {\em overconvergent} if for every $Y \to X$ in $X_{\acute{e}tale}$ with $Y$ affinoid the following holds: 
 
({\Large $\ast$}) \qquad 
For every special subset $V \subset Y$ we have \qquad 
$\DS \CS(V) \enspace = \enspace \lim_{\underset{V\subset\subset_Y U}{\longrightarrow}} \CS(U)\,,$ 
 
where the limit is taken over all special subsets $U \subset Y$ with $V\subset\subset_Y U$.  
\end{Def} 
 
The following is taken from \cite[{\S} 3]{JongPut}: An {\'e}tale sheaf 
$\CS$ is overconvergent if and only if for every $Y \to X$ in 
$X_{\acute{e}tale}$ there exists an admissible affinoid covering 
$Y_i,\,i\in I$, of $Y$ such that ({\Large $\ast$}) holds on each 
$Y_i$. 
The full subcategory 
of overconvergent {\'e}tale sheaves in the category of all abelian 
{\'e}tale sheaves is closed under the formation of quotients, subsheaves 
and extensions, i.e.,  it is a Serre subcategory. 
The space of {\em {\'e}tale analytic points} of $X$ is defined as  
\[ 
\CM_\et(X)\es:=\es\dirlim\Big\{\CM(Y):\es \text{for }Y\to X\hbox{ {\'e}tale}\Big\}\,. 
\] 
In line with \cite[{\S} 3.3]{JongPut}, {\'e}tale points $y$ can be assigned a
complete, algebraically closed residue field $k(y)$. We also regard $y$ as the
general morphism $y\!:\Spm k(y)\to X$.

For an {\'e}tale sheaf $\CS$ and $y\in\CM_\et(X)$ one may define a stalk 
$\CS_y$. {\'E}tale analytic points are sufficient to detect 
exactness in short exact sequences of overconvergent {\'e}tale 
sheaves. 
 
 
\subsection{Fundamental Groups} \label{SectFundamentalGroups} 
 
In this section we want to recall de Jong's \cite{deJong} definition of the {\em {\'e}tale fundamental group\/} and Andr{\'e}'s \cite[{\S} III.2]{andre2} definition of the {\em temperate {\'e}tale fundamental group} of a rigid analytic space. These concepts are best defined using Berkovich's analytic spaces; cf.\ \cite{berkovich}, \cite{berkovich2}. To stay within the framework of rigid analytic spaces, we therefore assume throughout this section that all rigid analytic spaces are quasi-separated and paracompact.

\medskip 
 
Let $X$ be an affinoid rigid analytic space over a complete, non-archimedean valued field $L$. 
 
\begin{Def} \label{DefFundamentalGroups} 
Let $\pi:Y\to X$ be a morphism of rigid analytic spaces over $L$. 
\begin{aufz} 
\item 
$Y$ is called an {\em {\'e}tale covering space of $X$\/} if every analytic point $x$ of $X$ has a wide affinoid neighborhood $U$ such that $Y\times_X U$ is a disjoint union of affinoids $V_i$ with $\pi|_{V_i}:V_i\to U$ finite {\'e}tale. 
\item 
$Y$ is called a {\em topological covering space of $X$\/} if we can choose $U$ and the $V_i$ above such that $\pi|_{V_i}:V_i\to U$ is an isomorphism. 
\item 
$Y$ is called a {\em finite {\'e}tale covering space of $X$\/} if $\pi$ is finite {\'e}tale 
\end{aufz} 
\end{Def}

The category of topological coverings is equivalent to the category of 
covering spaces of the topological Hausdorff space $\CM(X)$, cf.~
\cite[Lemma 2.6]{deJong}. The finite {\'e}tale covering spaces are 
used in algebraic geometry to define the algebraic fundamental group 
\cite{SGA1}.

\begin{Def} \label{DefTemperateFG} 
Let $Y$ be an {\'e}tale covering space of $X$. $Y$ is called a {\em temperate {\'e}tale covering space of $X$\/} if there exists a finite {\'e}tale covering space $S$ of $X$ and a topological covering space $T$ of $S$ such that $Y$ is the quotient of $T$ by an equivalence relation $R\subset T\times_X T$ which is a union of connected components. 
\end{Def} 
 
All the above notions of covering spaces are stable under taking 
connected components, fiber products and quotients.  
 
Clearly algebraic and topological coverings are temperate and 
temperate coverings are {\'e}tale. 
Furthermore, all {\'e}tale covering spaces are coverings for the {\'e}tale topology in the site $X_{\acute{e}tale}$; cf.\ Def.\ \ref{DefEtCov}. The reader should not confuse the concepts of {\'e}tale covering spaces and coverings for the {\'e}tale topology.

\end{appendix}

 
{

\vspace{2cm}

\parbox[t]{8.2cm}{ 
Gebhard B{\"o}ckle\\ 
Institut f{\"u}r Experimentelle Mathematik \\ 
Universit{\"a}t Duisburg-Essen, Campus Essen \\ 
Ellernstr.\ 29 \\ 
D -- 45326 Essen \\ 
Germany  \\[0.1cm] 
e-mail: boeckle@exp-math.uni-essen.de} 
\parbox[t]{6.7cm}{Urs Hartl  \\ 
Mathematisches Institut  \\ 
Albert-Ludwigs-Universit{\"a}t Freiburg  \\ 
Eckerstr.\ 1  \\ 
D -- 79104 Freiburg\\ 
Germany  \\[0.1cm] 
e-mail: urs.hartl@math.uni-freiburg.de 
} 
 
} 
 
\end{document}